# LOCAL RADEMACHER COMPLEXITIES

By Peter L. Bartlett, Olivier Bousquet and Shahar Mendelson

*University of California at Berkeley, Max Planck Institute for Biological Cybernetics and Australian National University*

We propose new bounds on the error of learning algorithms in terms of a data-dependent notion of complexity. The estimates we establish give optimal rates and are based on a local and empirical version of Rademacher averages, in the sense that the Rademacher averages are computed from the data, on a subset of functions with small empirical error. We present some applications to classification and prediction with convex function classes, and with kernel classes in particular.

**1. Introduction.** Estimating the performance of statistical procedures is useful for providing a better understanding of the factors that influence their behavior, as well as for suggesting ways to improve them. Although asymptotic analysis is a crucial first step toward understanding the behavior, finite sample error bounds are of more value as they allow the design of model selection (or parameter tuning) procedures. These error bounds typically have the following form: with high probability, the error of the estimator (typically a function in a certain class) is bounded by an empirical estimate of error plus a penalty term depending on the complexity of the class of functions that can be chosen by the algorithm. The differences between the true and empirical errors of functions in that class can be viewed as an empirical process. Many tools have been developed for understanding the behavior of such objects, and especially for evaluating their suprema—which can be thought of as a measure of how hard it is to estimate functions in the class at hand. The goal is thus to obtain the sharpest possible estimates on the complexity of function classes. A problem arises since the notion of complexity might depend on the (unknown) underlying probability measure according to which the data is produced. Distribution-free notions of the complexity, such as the Vapnik–Chervonenkis dimension [35] or the metric entropy [28], typically give conservative estimates. Distribution-dependent estimates, based for example on entropy numbers in the $L_2(P)$ distance, where $P$ is the underlying distribution, are not useful when $P$ is unknown. Thus, it is desirable to obtain data-dependent estimates which can readily be computed from the sample.

One of the most interesting data-dependent complexity estimates is the so-called Rademacher average associated with the class. Although known for a long time to be related to the expected supremum of the empirical process (thanks to symmetrization inequalities), it was first proposed as an effective complexity measure by Koltchinskii [15], Bartlett, Boucheron and Lugosi [1] and Mendelson [25] and then further studied in [3]. Unfortunately, one of the shortcomings of the Rademacher averages is that they provide *global* estimates of the complexity of the function class, that is, they do not reflect the fact that the algorithm will likely pick functions that have a small error, and in particular, only a small subset of the function class will be used. As a result, the best error rate that can be obtained via the global Rademacher averages is at least of the order of $1/\sqrt{n}$ (where $n$ is the sample size), which is suboptimal in some situations. Indeed, the type of algorithms we consider here are known in the statistical literature as $M$-estimators. They minimize an empirical loss criterion in a fixed class of functions. They have been extensively studied and their rate of convergence is known to be related to the modulus of continuity of the empirical process associated with the class of functions (rather than to the expected supremum of that empirical process). This modulus of continuity is well understood from the empirical processes theory viewpoint (see, e.g., [33, 34]). Also, from the point of view of $M$-estimators, the quantity which determines the rate of convergence is actually a fixed point of this modulus of continuity. Results of this type have been obtained by van de Geer [31, 32] (among others), who also provides nonasymptotic exponential inequalities. Unfortunately, these are in terms of entropy (or random entropy) and hence might not be useful when the probability distribution is unknown.

The key property that allows one to prove fast rates of convergence is the fact that around the best function in the class, the variance of the increments of the empirical process [or the $L_2(P)$ distance to the best function] is upper bounded by a linear function of the expectation of these increments. In the context of regression with squared loss, this happens as soon as the functions are bounded and the class of functions is convex. In the context of classification, Mammen and Tsybakov have shown [20] that this also happens under conditions on the conditional distribution (especially about its behavior around $1/2$). They actually do not require the relationship between variance and expectation (of the increments) to be linear but allow for more general, power type inequalities. Their results, like those of van de Geer, are asymptotic.

In order to exploit this key property and have finite sample bounds, rather than considering the Rademacher averages of the entire class as the complexity measure, it is possible to consider the Rademacher averages of a small subset of the class, usually the intersection of the class with a ball centered at a function of interest. These *local* Rademacher averages can serve as a complexity measure; clearly, they are always smaller than the corresponding global averages. Several authors have considered the use of local estimates of the complexity of the function class in order to obtain better bounds. Before presenting their results, we introduce some notation which is used throughout the paper.

Let $(\mathcal{X}, P)$ be a probability space. Denote by $\mathcal{F}$ a class of measurable functions from $\mathcal{X}$ to $\mathbb{R}$, and set $X_1, \ldots, X_n$ to be independent random variables distributed according to $P$. Let $\sigma_1, \ldots, \sigma_n$ be $n$ independent *Rademacher* random variables, that is, independent random variables for which $\Pr(\sigma_i = 1) = \Pr(\sigma_i = -1) = 1/2$.

For a function $f : \mathcal{X} \to \mathbb{R}$, define

$$P_n f = \frac{1}{n} \sum_{i=1}^n f(X_i), \qquad Pf = \mathbb{E} f(X), \qquad R_n f = \frac{1}{n} \sum_{i=1}^n \sigma_i f(X_i).$$

For a class $\mathcal{F}$, set

$$R_n \mathcal{F} = \sup_{f \in \mathcal{F}} R_n f.$$

Define $\mathbb{E}_\sigma$ to be the expectation with respect to the random variables $\sigma_1, \ldots, \sigma_n$, conditioned on all of the other random variables. The Rademacher average of $\mathcal{F}$ is $\mathbb{E} R_n \mathcal{F}$, and the empirical (or conditional) Rademacher averages of $\mathcal{F}$ are

$$\mathbb{E}_\sigma R_n \mathcal{F} = \frac{1}{n} \mathbb{E} \left( \sup_{f \in \mathcal{F}} \sum_{i=1}^n \sigma_i f(X_i) | X_1, \ldots, X_n \right).$$

Some classical properties of Rademacher averages and some simple lemmas (which we use often) are listed in Appendix A.1.

The simplest way to obtain the property allowing for fast rates of convergence is to consider nonnegative uniformly bounded functions (or increments with respect to a fixed null function). In this case, one trivially has for all $f \in \mathcal{F}$, $\mathsf{Var}[f] \leq cPf$. This is exploited by Koltchinskii and Panchenko [16], who consider the case of prediction with absolute loss when functions in $\mathcal{F}$ have values in $[0, 1]$ and there is a *perfect* function $f^*$ in the class, that is, $Pf^* = 0$. They introduce an iterative method involving local empirical Rademacher averages. They first construct a function $\phi_n(r) = c_1 R_n \{f : P_n f \leq 2r\} + c_2 \sqrt{rx/n} +$ $c_3/n$, which can be computed from the data. For $\hat{r}_N$ defined by $\hat{r}_0 = 1$, $\hat{r}_{k+1} = \phi_n(\hat{r}_k)$, they show that with probability at least $1 - 2Ne^{-x}$,

$$P\hat{f} \leq \hat{r}_N + \frac{2x}{n},$$

where $\hat{f}$ is a minimizer of the empirical error, that is, a function in $\mathcal{F}$ satisfying $P_n\hat{f} = \inf_{f \in \mathcal{F}} P_n f$. Hence, this nonincreasing sequence of local Rademacher averages can be used as upper bounds on the error of the empirical minimizer $\hat{f}$. Furthermore, if $\psi_n$ is a concave function such that $\psi_(\sqrt{r}) \geq \mathbb{E}_\sigma R_n\{f \in \mathcal{F} : P_n f \leq r\}$, and if the number of iterations $N$ is at least $1 + \lceil \log_2 \log_2 n/x \rceil$, then with probability at least $1 - Ne^{-x}$,

$$\hat{r}_N \leq c\left(\hat{r}^* + \frac{x}{n}\right),$$

where $r^*$ is a solution of the fixed-point equation $\psi(\sqrt{r}) = r$. Combining the above results, one has a procedure to obtain data-dependent error bounds that are of the order of the fixed point of the modulus of continuity at 0 of the empirical Rademacher averages. One limitation of this result is that it assumes that there is a function $f^*$ in the class with $Pf^* = 0$. In contrast, we are interested in prediction problems where $Pf$ is the error of an estimator, and in the presence of noise there may not be any perfect estimator (even the best in the class can have nonzero error).

More recently, Bousquet, Koltchinskii and Panchenko [9] have obtained a more general result avoiding the iterative procedure. Their result is that for functions with values in $[0,1]$, with probability at least $1 - e^{-x}$,

$$(1.1) \qquad \forall f \in \mathcal{F} \qquad Pf \leq c\left(P_n f + \hat{r}^* + \frac{t + \log \log n}{n}\right),$$

where $\hat{r}^*$ is the fixed point of a concave function $\psi_n$ satisfying $\psi_n(0) = 0$ and

$$\psi_n(\sqrt{r}) \geq \mathbb{E}_\sigma R_n\{f \in \mathcal{F} : P_n f \leq r\}.$$

The main difference between this and the results of [16] is that there is no requirement that the class contain a perfect function. However, the local Rademacher averages are centered around the zero function instead of the one that minimizes $Pf$. As a consequence, the fixed point $\hat{r}^*$ cannot be expected to converge to zero when $\inf_{f \in \mathcal{F}} Pf > 0$.

In order to remove this limitation, Lugosi and Wegkamp [19] use localized Rademacher averages of a small ball around the minimizer $\hat{f}$ of $P_n$. However, their result is restricted to nonnegative functions, and in particular functions with values in $\{0,1\}$. Moreover, their bounds also involve some global information, in the form of the shatter coefficients $S_\mathcal{F}(X_1^n)$ of the function class (i.e., the cardinality of the *coordinate projections* of the class $\mathcal{F}$ on the data $X_1^n$). They show that there are constants $c_1, c_2$ such that, with probability at least $1 - 8/n$, the empirical minimizer $\hat{f}$ satisfies

$$P\hat{f} \leq \inf_{f \in \mathcal{F}} Pf + 2\widehat{\psi}_n(\hat{r}_n),$$

where

$$\widehat{\psi}_n(r) = c_1 \left( \mathbb{E}_\sigma R_n \{f \in \mathcal{F} : P_n f \leq 16 P_n \hat{f} + 15r\} + \frac{\log n}{n} + \sqrt{\frac{\log n}{n}} \sqrt{P_n \hat{f} + r} \right)$$

and $\hat{r}_n = c_2(\log S_\mathcal{F}(X_1^n) + \log n)/n$. The limitation of this result is that $\hat{r}_n$ has to be chosen according to the (empirically measured) complexity of the whole class, which may not be as sharp as the Rademacher averages, and in general, is not a fixed point of $\widehat{\psi}_n$. Moreover, the balls over which the Rademacher averages are computed in $\widehat{\psi}_n$ contain a factor of 16 in front of $P_n \hat{f}$. As we explain later, this induces a lower bound on $\widehat{\psi}_n$ when there is no function with $Pf = 0$ in the class.

It seems that the only way to capture the right behavior in the general, noisy case is to analyze the increments of the empirical process, in other words, to directly consider the functions $f - f^*$. This approach was first proposed by Massart [22]; see also [26]. Massart introduces the assumption

$$\mathsf{Var}[\ell_f(X) - \ell_{f^*}(X)] \leq d^2(f, f^*) \leq B(P\ell_f - P\ell_{f^*}),$$

where $\ell_f$ is the loss associated with the function $f$ [in other words, $\ell_f(X, Y) = \ell(f(X), Y)$, which measures the discrepancy in the prediction made by $f$], $d$ is a pseudometric and $f^*$ minimizes the expected loss. (The previous results could also be stated in terms of loss functions, but we omitted this in order to simplify exposition. However, the extra notation is necessary to properly state Massart's result.) This is a more refined version of the assumption we mentioned earlier on the relationship between the variance and expectation of the increments of the empirical process. It is only satisfied for some loss functions $\ell$ and function classes $\mathcal{F}$. Under this assumption, Massart considers a nondecreasing function $\psi$ satisfying

$$\psi(r) \geq \mathbb{E} \sup_{f \in \mathcal{F}, d^2(f, f^*)^2 \leq r} |Pf - Pf^* - P_n f + P_n f^*| + c \frac{x}{n},$$

such that $\psi(r)/\sqrt{r}$ is nonincreasing (we refer to this property as the sub-root property later in the paper). Then, with probability at least $1 - e^{-x}$,

(1.2) $$\forall f \in \mathcal{F} \qquad P\ell_f - P\ell_{f^*} \leq c\left(r^* + \frac{x}{n}\right),$$

where $r^*$ is the fixed point of $\psi$ and $c$ depends only on $B$ and on the uniform bound on the range of functions in $\mathcal{F}$. It can be proved that in many situations of interest, this bound suffices to prove minimax rates of convergence for penalized $M$-estimators. (Massart considers examples where the complexity term can be bounded using a priori global information about the function class.) However, the main limitation of this result is that it does not involve quantities that can be computed from the data.

Finally, as we mentioned earlier, Mendelson [26] gives an analysis similar to that of Massart, in a slightly less general case (with no noise in the target values, i.e., the conditional distribution of $Y$ given $X$ is concentrated at one point). Mendelson introduces the notion of the *star-hull* of a class of functions (see the next section for a definition) and considers Rademacher averages of this star-hull as a localized measure of complexity. His results also involve a priori knowledge of the class, such as the rate of growth of covering numbers.

We can now spell out our goal in more detail: in this paper we combine the increment-based approach of Massart and Mendelson (dealing with differences of functions, or more generally with bounded real-valued functions) with the empirical local Rademacher approach of Koltchinskii and Panchenko and of Lugosi and Wegkamp, in order to obtain data-dependent bounds which depend on a fixed point of the modulus of continuity of Rademacher averages computed around the empirically best function.

Our first main result (Theorem 3.3) is a distribution-dependent result involving the fixed point $r^*$ of a local Rademacher average of the star-hull of the class $\mathcal{F}$. This shows that functions with the sub-root property can readily be obtained from Rademacher averages, while in previous work the appropriate functions were obtained only via global information about the class.

The second main result (Theorems 4.1 and 4.2) is an empirical counterpart of the first one, where the complexity is the fixed point of an empirical local Rademacher average. We also show that this fixed point is within a constant factor of the nonempirical one.

Equipped with this result, we can then prove (Theorem 5.4) a fully data-dependent analogue of Massart's result, where the Rademacher averages are localized around the minimizer of the empirical loss.

We also show (Theorem 6.3) that in the context of classification, the local Rademacher averages of star-hulls can be approximated by solving a weighted empirical error minimization problem.

Our final result (Corollary 6.7) concerns regression with kernel classes, that is, classes of functions that are generated by a positive definite kernel. These classes are widely used in interpolation and estimation problems as they yield computationally efficient algorithms. Our result gives a data-dependent complexity term that can be computed directly from the eigenvalues of the Gram matrix (the matrix whose entries are values of the kernel on the data).

The sharpness of our results is demonstrated from the fact that we recover, in the distribution-dependent case (treated in Section 4), similar results to those of Massart [22], which, in the situations where they apply, give the minimax optimal rates or the best known results. Moreover, the data-dependent bounds that we obtain as counterparts of these results have the same rate of convergence (see Theorem 4.2).

The paper is organized as follows. In Section 2 we present some preliminary results obtained from concentration inequalities, which we use throughout. Section 3 establishes error bounds using local Rademacher averages and explains how to compute their fixed points from "global information" (e.g., estimates of the metric entropy or of the combinatorial dimensions of the indexing class), in which case the optimal estimates can be recovered. In Section 4 we give a data-dependent error bound using empirical and local Rademacher averages, and show the connection between the fixed points of the empirical and nonempirical Rademacher averages. In Section 5 we apply our results to loss classes. We give estimates that generalize the results of Koltchinskii and Panchenko by eliminating the requirement that some function in the class have zero loss, and are more general than those of Lugosi and Wegkamp, since there is no need have in our case to estimate global shatter coefficients of the class. We also give a data-dependent extension of Massart's result where the local averages are computed around the minimizer of the empirical loss. Finally, Section 6 shows that the problem of estimating these local Rademacher averages in classification reduces to weighted empirical risk minimization. It also shows that the local averages for kernel classes can be sharply bounded in terms of the eigenvalues of the Gram matrix.

**2. Preliminary results.** Recall that the star-hull of $\mathcal{F}$ around $f_0$ is defined by

$$\mathrm{star}(\mathcal{F}, f_0) = \{f_0 + \alpha(f - f_0) : f \in \mathcal{F}, \alpha \in [0,1]\}.$$

Throughout this paper, we will manipulate suprema of empirical processes, that is, quantities of the form $\sup_{f \in \mathcal{F}}(Pf - P_n f)$. We will always assume they are measurable without explicitly mentioning it. In other words, we assume that the class $\mathcal{F}$ and the distribution $P$ satisfy appropriate (mild) conditions for measurability of this supremum (we refer to [11, 28] for a detailed account of such issues).

The following theorem is the main result of this section and is at the core of all the proofs presented later. It shows that if the functions in a class have small variance, the maximal deviation between empirical means and true means is controlled by the Rademacher averages of $\mathcal{F}$. In particular, the bound improves as the largest variance of a class member decreases.

THEOREM 2.1. *Let $\mathcal{F}$ be a class of functions that map $\mathcal{X}$ into $[a,b]$. Assume that there is some $r > 0$ such that for every $f \in \mathcal{F}$, $\mathsf{Var}[f(X_i)] \leq r$. Then, for every $x > 0$, with probability at least $1 - e^{-x}$,*

$$\sup_{f \in \mathcal{F}}(Pf - P_n f) \leq \inf_{\alpha > 0}\left(2(1+\alpha)\mathbb{E}R_n\mathcal{F} + \sqrt{\frac{2rx}{n}} + (b-a)\left(\frac{1}{3} + \frac{1}{\alpha}\right)\frac{x}{n}\right),$$

*and with probability at least $1 - 2e^{-x}$,*

$$\sup_{f \in \mathcal{F}}(Pf - P_n f)$$
$$\leq \inf_{\alpha \in (0,1)}\left(2\frac{1+\alpha}{1-\alpha}\mathbb{E}_\sigma R_n\mathcal{F} + \sqrt{\frac{2rx}{n}} + (b-a)\left(\frac{1}{3} + \frac{1}{\alpha} + \frac{1+\alpha}{2\alpha(1-\alpha)}\right)\frac{x}{n}\right).$$

*Moreover, the same results hold for the quantity $\sup_{f \in \mathcal{F}}(P_n f - Pf)$.*

This theorem, which is proved in Appendix A.2, is a more or less direct consequence of Talagrand's inequality for empirical processes [30]. However, the actual statement presented here is new in the sense that it displays the best known constants. Indeed, compared to the previous result of Koltchinskii and Panchenko [16] which was based on Massart's version of Talagrand's inequality [21], we have used the most refined concentration inequalities available: that of Bousquet [7] for the supremum of the empirical process and that of Boucheron, Lugosi and Massart [5] for the Rademacher averages. This last inequality is a powerful tool to obtain data-dependent bounds, since it allows one to replace the Rademacher average (which measures the complexity of the class of functions) by its empirical version, which can be efficiently computed in some cases. Details about these inequalities are given in Appendix A.1.

When applied to the full function class $\mathcal{F}$, the above theorem is not useful. Indeed, with only a trivial bound on the maximal variance, better results can be obtained via simpler concentration inequalities, such as the bounded difference inequality [23], which would allow $\sqrt{rx/n}$ to be replaced by $\sqrt{x/n}$. However, by applying Theorem 2.1 to subsets of $\mathcal{F}$ or to modified classes obtained from $\mathcal{F}$, much better results can be obtained. Hence, the presence of an upper bound on the variance in the square root term is the key ingredient of this result.

A last preliminary result that we will require is the following consequence of Theorem 2.1, which shows that if the local Rademacher averages are small, then balls in $L_2(P)$ are probably contained in the corresponding empirical balls [i.e., in $L_2(P_n)$] with a slightly larger radius.

COROLLARY 2.2. *Let $\mathcal{F}$ be a class of functions that map $\mathcal{X}$ into $[-b,b]$ with $b > 0$. For every $x > 0$ and $r$ that satisfy*

$$r \geq 10b\mathbb{E}R_n\{f : f \in \mathcal{F}, Pf^2 \leq r\} + \frac{11b^2 x}{n},$$

*then with probability at least* $1 - e^{-x}$,

$$\{f \in \mathcal{F} : Pf^2 \leq r\} \subseteq \{f \in \mathcal{F} : P_n f^2 \leq 2r\}.$$

PROOF. Since the range of any function in the set $\mathcal{F}_r = \{f^2 : f \in \mathcal{F}, Pf^2 \leq r\}$ is contained in $[0, b^2]$, it follows that $\mathsf{Var}[f^2(X_i)] \leq Pf^4 \leq b^2 Pf^2 \leq b^2 r$. Thus, by the first part of Theorem 2.1 (with $\alpha = 1/4$), with probability at least $1 - e^{-x}$, every $f \in \mathcal{F}_r$ satisfies

$$P_n f^2 \leq r + \frac{5}{2} \mathbb{E} R_n \{f^2 : f \in \mathcal{F}, Pf^2 \leq r\} + \sqrt{\frac{2b^2 r x}{n}} + \frac{13 b^2 x}{3n}$$

$$\leq r + \frac{5}{2} \mathbb{E} R_n \{f^2 : f \in \mathcal{F}, Pf^2 \leq r\} + \frac{r}{2} + \frac{16 b^2 x}{3n}$$

$$\leq r + 5b \mathbb{E} R_n \{f : f \in \mathcal{F}, Pf^2 \leq r\} + \frac{r}{2} + \frac{16 b^2 x}{3n}$$

$$\leq 2r,$$

where the second inequality follows from Lemma A.3 and we have used, in the second to last inequality, Theorem A.6 applied to $\phi(x) = x^2$ (with Lipschitz constant $2b$ on $[-b, b]$). □

**3. Error bounds with local complexity.** In this section we show that the Rademacher averages associated with a small subset of the class may be considered as a complexity term in an error bound. Since these *local* Rademacher averages are always smaller than the corresponding global averages, they lead to sharper bounds.

We present a general error bound involving local complexities that is applicable to classes of bounded functions for which the variance is bounded by a fixed linear function of the expectation. In this case the local Rademacher averages are defined as $\mathbb{E} R_n \{f \in \mathcal{F} : T(f) \leq r\}$ where $T(f)$ is an upper bound on the variance [typically chosen as $T(f) = Pf^2$].

There is a trade-off between the size of the subset we consider in these local averages and its complexity; we shall see that the optimal choice is given by a fixed point of an upper bound on the local Rademacher averages. The functions we use as upper bounds are *sub-root* functions; among other useful properties, sub-root functions have a unique fixed point.

DEFINITION 3.1. *A function* $\psi : [0, \infty) \to [0, \infty)$ *is* sub-root *if it is nonnegative, nondecreasing and if* $r \mapsto \psi(r)/\sqrt{r}$ *is nonincreasing for* $r > 0$.

We only consider nontrivial sub-root functions, that is, sub-root functions that are not the constant function $\psi \equiv 0$.

LEMMA 3.2. *If $\psi:[0,\infty) \to [0,\infty)$ is a nontrivial sub-root function, then it is continuous on $[0,\infty)$ and the equation $\psi(r) = r$ has a unique positive solution. Moreover, if we denote the solution by $r^*$, then for all $r > 0$, $r \geq \psi(r)$ if and only if $r^* \leq r$.*

The proof of this lemma is in Appendix A.2. In view of the lemma, we will simply refer to the quantity $r^*$ as the *unique positive solution* of $\psi(r) = r$, or as the *fixed point* of $\psi$.

3.1. *Error bounds.* We can now state and discuss the main result of this section. It is composed of two parts: in the first part, one requires a sub-root upper bound on the local Rademacher averages, and in the second part, it is shown that better results can be obtained when the class over which the averages are computed is enlarged slightly.

THEOREM 3.3. *Let $\mathcal{F}$ be a class of functions with ranges in $[a,b]$ and assume that there are some functional $T:\mathcal{F} \to \mathbb{R}^+$ and some constant $B$ such that for every $f \in \mathcal{F}$, $\mathsf{Var}[f] \leq T(f) \leq BPf$. Let $\psi$ be a sub-root function and let $r^*$ be the fixed point of $\psi$.*

1. *Assume that $\psi$ satisfies, for any $r \geq r^*$,*

$$\psi(r) \geq B\mathbb{E}R_n\{f \in \mathcal{F} : T(f) \leq r\}.$$

*Then, with $c_1 = 704$ and $c_2 = 26$, for any $K > 1$ and every $x > 0$, with probability at least $1 - e^{-x}$,*

$$\forall f \in \mathcal{F} \qquad Pf \leq \frac{K}{K-1}P_n f + \frac{c_1 K}{B}r^* + \frac{x(11(b-a) + c_2 BK)}{n}.$$

*Also, with probability at least $1 - e^{-x}$,*

$$\forall f \in \mathcal{F} \qquad P_n f \leq \frac{K+1}{K}Pf + \frac{c_1 K}{B}r^* + \frac{x(11(b-a) + c_2 BK)}{n}.$$

2. *If, in addition, for $f \in \mathcal{F}$ and $\alpha \in [0,1]$, $T(\alpha f) \leq \alpha^2 T(f)$, and if $\psi$ satisfies, for any $r \geq r^*$,*

$$\psi(r) \geq B\mathbb{E}R_n\{f \in \mathrm{star}(\mathcal{F},0) : T(f) \leq r\},$$

*then the same results hold true with $c_1 = 6$ and $c_2 = 5$.*

The proof of this theorem is given in Section 3.2.

We can compare the results to our starting point (Theorem 2.1). The improvement comes from the fact that the complexity term, which was essentially $\sup_r \psi(r)$ in Theorem 2.1 (if we had applied it to the class $\mathcal{F}$ directly) is now reduced to $r^*$, the fixed point of $\psi$. So the complexity term is always smaller (later, we show how to estimate $r^*$). On the other hand, there is some loss since the constant in front of $P_n f$ is strictly larger than 1. Section 5.2 will show that this is not an issue in the applications we have in mind.

In Sections 5.1 and 5.2 we investigate conditions that ensure the assumptions of this theorem are satisfied, and we provide applications of this result to prediction problems. The condition that the variance is upper bounded by the expectation turns out to be crucial to obtain these results.

The idea behind Theorem 3.3 originates in the work of Massart [22], who proves a slightly different version of the first part. The difference is that we use local Rademacher averages instead of the expectation of the supremum of the empirical process on a ball. Moreover, we give smaller constants. As far as we know, the second part of Theorem 3.3 is new.

3.1.1. *Choosing the function $\psi$.* Notice that the function $\psi$ cannot be chosen arbitrarily and has to satisfy the sub-root property. One possible approach is to use classical upper bounds on the Rademacher averages, such as Dudley's entropy integral. This can give a sub-root upper bound and was used, for example, in [16] and in [22].

However, the second part of Theorem 3.3 indicates a possible choice for $\psi$, namely, one can take $\psi$ as the local Rademacher averages of the star-hull of $\mathcal{F}$ around 0. The reason for this comes from the following lemma, which shows that if the class is star-shaped and $T(f)$ behaves as a quadratic function, the Rademacher averages are sub-root.

LEMMA 3.4. *If the class $\mathcal{F}$ is star-shaped around $\hat{f}$ (which may depend on the data), and $T : \mathcal{F} \to \mathbb{R}^+$ is a (possibly random) function that satisfies $T(\alpha f) \leq \alpha^2 T(f)$ for any $f \in \mathcal{F}$ and any $\alpha \in [0, 1]$, then the (random) function $\psi$ defined for $r \geq 0$ by*

$$\psi(r) = \mathbb{E}_\sigma R_n \{ f \in \mathcal{F} : T(f - \hat{f}) \leq r \}$$

*is sub-root and $r \mapsto \mathbb{E}\psi(r)$ is also sub-root.*

This lemma is proved in Appendix A.2.

Notice that making a class star-shaped only increases it, so that

$$\mathbb{E} R_n \{ f \in \mathrm{star}(\mathcal{F}, f_0) : T(f) \leq r \} \geq \mathbb{E} R_n \{ f \in \mathcal{F} : T(f) \leq r \}.$$

However, this increase in size is moderate as can be seen, for example, if one compares covering numbers of a class and its star-hull (see, e.g., [26], Lemma 4.5).

3.1.2. *Some consequences.* As a consequence of Theorem 3.3, we obtain an error bound when $\mathcal{F}$ consists of uniformly bounded nonnegative functions. Notice that in this case the variance is trivially bounded by a constant times the expectation and one can directly use $T(f) = Pf$.

COROLLARY 3.5. *Let $\mathcal{F}$ be a class of functions with ranges in $[0, 1]$. Let $\psi$ be a sub-root function, such that for all $r \geq 0$,*

$$\mathbb{E}R_n\{f \in \mathcal{F} : Pf \leq r\} \leq \psi(r),$$

*and let $r^*$ be the fixed point of $\psi$. Then, for any $K > 1$ and every $x > 0$, with probability at least $1 - e^{-x}$, every $f \in \mathcal{F}$ satisfies*

$$Pf \leq \frac{K}{K-1} P_n f + 704 K r^* + \frac{x(11 + 26K)}{n}.$$

*Also, with probability at least $1 - e^{-x}$, every $f \in \mathcal{F}$ satisfies*

$$P_n f \leq \frac{K+1}{K} P f + 704 K r^* + \frac{x(11 + 26K)}{n}.$$

PROOF. When $f \in [0, 1]$, we have $\mathsf{Var}[f] \leq Pf$ so that the result follows from applying Theorem 3.3 with $T(f) = Pf$. □

We also note that the same idea as in the proof of Theorem 3.3 gives a converse of Corollary 2.2, namely, that with high probability the intersection of $\mathcal{F}$ with an empirical ball of a fixed radius is contained in the intersection of $\mathcal{F}$ with an $L_2(P)$ ball with a slightly larger radius.

LEMMA 3.6. *Let $\mathcal{F}$ be a class of functions that map $\mathcal{X}$ into $[-1, 1]$. Fix $x > 0$. If*

$$r \geq 20 \mathbb{E} R_n \{f : f \in \text{star}(\mathcal{F}, 0), Pf^2 \leq r\} + \frac{26x}{n},$$

*then with probability at least $1 - e^{-x}$,*

$$\{f \in \text{star}(\mathcal{F}, 0) : P_n f^2 \leq r\} \subseteq \{f \in \text{star}(\mathcal{F}, 0) : Pf^2 \leq 2r\}.$$

This result, proved in Section 3.2, will be useful in Section 4.

3.1.3. *Estimating $r^*$ from global information.* The error bounds involve fixed points of functions that define upper bounds on the local Rademacher averages. In some cases these fixed points can be estimated from global information on the function class. We present a complete analysis only in a simple case, where $\mathcal{F}$ is a class of binary-valued functions with a finite VC-dimension.

COROLLARY 3.7. *Let $\mathcal{F}$ be a class of $\{0,1\}$-valued functions with VC-dimension $d < \infty$. Then for all $K > 1$ and every $x > 0$, with probability at least $1 - e^{-x}$, every $f \in \mathcal{F}$ satisfies*

$$Pf \leq \frac{K}{K-1} P_n f + cK \left( \frac{d \log(n/d)}{n} + \frac{x}{n} \right).$$

The proof is in Appendix A.2.

The above result is similar to results obtained by Vapnik and Chervonenkis [35] and by Lugosi and Wegkamp (Theorem 3.1 of [19]). However, they used inequalities for weighted empirical processes indexed by nonnegative functions. Our results have more flexibility since they can accommodate general functions, although this is not needed in this simple corollary.

The proof uses a similar line of reasoning to proofs in [26, 27]. Clearly, it extends to any class of real-valued functions for which one has estimates for the entropy integral, such as classes with finite pseudo-dimension or a combinatorial dimension that grows more slowly than quadratically. See [26, 27] for more details.

Notice also that the rate of $\log n / n$ is the best known.

3.1.4. *Proof techniques.* Before giving the proofs of the results mentioned above, let us sketch the techniques we use. The approach has its roots in classical empirical processes theory, where it was understood that the modulus of continuity of the empirical process is an important quantity (here $\psi$ plays this role). In order to obtain nonasymptotic results, two approaches have been developed: the first one consists of cutting the class $\mathcal{F}$ into smaller pieces, where one has control of the variance of the elements. This is the so-called *peeling* technique (see, e.g., [31, 32, 33, 34] and references therein). The second approach consists of weighting the functions in $\mathcal{F}$ by dividing them by their variance. Many results have been obtained on such weighted empirical processes (see, e.g., [28]). The results of Vapnik and Chervonenkis based on weighting [35] are restricted to classes of nonnegative functions. Also, most previous results, such as those of Pollard [28], van de Geer [32] or Haussler [13], give complexity terms that involve "global" measures of complexity of the class, such as covering numbers. None of these results uses the recently introduced Rademacher averages as measures of complexity. It turns out that it is possible to combine the peeling and weighting ideas with concentration inequalities to obtain such results, as proposed by Massart in [22], and also used (for nonnegative functions) by Koltchinskii and Panchenko [16].

The idea is the following:

(a) Apply Theorem 2.1 to the class of functions $\{f/w(f) : f \in \mathcal{F}\}$, where $w$ is some nonnegative weight of the order of the variance of $f$. Hence, the functions in this class have a small variance.

(b) Upper bound the Rademacher averages of this weighted class, by "peeling off" subclasses of $\mathcal{F}$ according to the variance of their elements, and bounding the Rademacher averages of these subclasses using $\psi$.

(c) Use the sub-root property of $\psi$, so that its fixed point gives a common upper bound on the complexity of all the subclasses (up to some scaling).

(d) Finally, convert the upper bound for functions in the weighted class into a bound for functions in the initial class.

The idea of peeling—that is, of partitioning the class $\mathcal{F}$ into slices where functions have variance within a certain range—is at the core of the proof of the first part of Theorem 3.3 [see, e.g., (3.1)]. However, it does not appear explicitly in the proof of the second part. One explanation is that when one considers the star-hull of the class, it is enough to consider two subclasses: the functions with $T(f) \leq r$ and the ones with $T(f) > r$, and this is done by introducing the weighting factor $T(f) \vee r$. This idea was exploited in the work of Mendelson [26] and, more recently, in [4]. Moreover, when one considers the set $\mathcal{F}_r = \text{star}(\mathcal{F}, 0) \cap \{T(f) \leq r\}$, any function $f' \in \mathcal{F}$ with $T(f') > r$ will have a scaled down representative in that set. So even though it seems that we look at the class $\text{star}(\mathcal{F}, 0)$ only locally, we still take into account all of the functions in $\mathcal{F}$ (with appropriate scaling).

3.2. *Proofs.* Before presenting the proof, let us first introduce some additional notation. Given a class $\mathcal{F}$, $\lambda > 1$ and $r > 0$, let $w(f) = \min\{r\lambda^k : k \in \mathbb{N}, r\lambda^k \geq T(f)\}$ and set

$$\mathcal{G}_r = \left\{ \frac{r}{w(f)} f : f \in \mathcal{F} \right\}.$$

Notice that $w(f) \geq r$, so that $\mathcal{G}_r \subseteq \{\alpha f : f \in \mathcal{F}, \alpha \in [0,1]\} = \text{star}(\mathcal{F}, 0)$. Define

$$V_r^+ = \sup_{g \in \mathcal{G}_r} Pg - P_n g \quad \text{and} \quad V_r^- = \sup_{g \in \mathcal{G}_r} P_n g - Pg.$$

For the second part of the theorem, we need to introduce another class of functions,

$$\tilde{\mathcal{G}}_r := \left\{ \frac{rf}{T(f) \vee r} : f \in \mathcal{F} \right\},$$

and define

$$\tilde{V}_r^+ = \sup_{g \in \tilde{\mathcal{G}}_r} Pg - P_n g \quad \text{and} \quad \tilde{V}_r^- = \sup_{g \in \tilde{\mathcal{G}}_r} P_n g - Pg.$$

LEMMA 3.8. *With the above notation, assume that there is a constant $B > 0$ such that for every $f \in \mathcal{F}$, $T(f) \leq BPf$. Fix $K > 1$, $\lambda > 0$ and $r > 0$.*

If $V_r^+ \leq r/(\lambda BK)$, then

$$\forall f \in \mathcal{F} \qquad Pf \leq \frac{K}{K-1} P_n f + \frac{r}{\lambda BK}.$$

Also, if $V_r^- \leq r/(\lambda BK)$, then

$$\forall f \in \mathcal{F} \qquad P_n f \leq \frac{K+1}{K} Pf + \frac{r}{\lambda BK}.$$

Similarly, if $K > 1$ and $r > 0$ are such that $\tilde{V}_r^+ \leq r/(BK)$, then

$$\forall f \in \mathcal{F} \qquad Pf \leq \frac{K}{K-1} P_n f + \frac{r}{BK}.$$

Also, if $\tilde{V}_r^- \leq r/(BK)$, then

$$\forall f \in \mathcal{F} \qquad P_n f \leq \frac{K+1}{K} Pf + \frac{r}{BK}.$$

PROOF. Notice that for all $g \in \mathcal{G}_r$, $Pg \leq P_n g + V_r^+$. Fix $f \in \mathcal{F}$ and define $g = rf/w(f)$. When $T(f) \leq r$, $w(f) = r$, so that $g = f$. Thus, the fact that $Pg \leq P_n g + V_r^+$ implies that $Pf \leq P_n f + V_r^+ \leq P_n f + r/(\lambda BK)$.

On the other hand, if $T(f) > r$, then $w(f) = r\lambda^k$ with $k > 0$ and $T(f) \in (r\lambda^{k-1}, r\lambda^k]$. Moreover, $g = f/\lambda^k$, $Pg \leq P_n g + V_r^+$, and thus

$$\frac{Pf}{\lambda^k} \leq \frac{P_n f}{\lambda^k} + V_r^+.$$

Using the fact that $T(f) > r\lambda^{k-1}$, it follows that

$$Pf \leq P_n f + \lambda^k V_r^+ < P_n f + \lambda T(f) V_r^+/r \leq P_n f + Pf/K.$$

Rearranging,

$$Pf \leq \frac{K}{K-1} P_n f < \frac{K}{K-1} P_n f + \frac{r}{\lambda BK}.$$

The proof of the second result is similar. For the third and fourth results, the reasoning is the same. □

PROOF OF THEOREM 3.3, FIRST PART. Let $\mathcal{G}_r$ be defined as above, where $r$ is chosen such that $r \geq r^*$, and note that functions in $\mathcal{G}_r$ satisfy $\|g - Pg\|_\infty \leq b - a$ since $0 \leq r/w(f) \leq 1$. Also, we have $\mathsf{Var}[g] \leq r$. Indeed, if $T(f) \leq r$, then $g = f$, and thus $\mathsf{Var}[g] = \mathsf{Var}[f] \leq r$. Otherwise, when $T(f) > r$, $g = f/\lambda^k$ (where $k$ is such that $T(f) \in (r\lambda^{k-1}, r\lambda^k]$), so that $\mathsf{Var}[g] = \mathsf{Var}[f]/\lambda^{2k} \leq r$.

Applying Theorem 2.1, for all $x > 0$, with probability $1 - e^{-x}$,

$$V_r^+ \leq 2(1+\alpha) \mathbb{E} R_n \mathcal{G}_r + \sqrt{\frac{2rx}{n}} + (b-a)\left(\frac{1}{3} + \frac{1}{\alpha}\right)\frac{x}{n}.$$

Let $\mathcal{F}(x,y) := \{f \in \mathcal{F} : x \leq T(f) \leq y\}$ and define $k$ to be the smallest integer such that $r\lambda^{k+1} \geq Bb$. Then

$$\mathbb{E}R_n\mathcal{G}_r \leq \mathbb{E}R_n\mathcal{F}(0,r) + \mathbb{E} \sup_{f \in \mathcal{F}(r,Bb)} \frac{r}{w(f)} R_n f$$

$$\leq \mathbb{E}R_n\mathcal{F}(0,r) + \sum_{j=0}^{k} \mathbb{E} \sup_{f \in \mathcal{F}(r\lambda^j, r\lambda^{j+1})} \frac{r}{w(f)} R_n f$$

(3.1)

$$= \mathbb{E}R_n\mathcal{F}(0,r) + \sum_{j=0}^{k} \lambda^{-j} \mathbb{E} \sup_{f \in \mathcal{F}(r\lambda^j, r\lambda^{j+1})} R_n f$$

$$\leq \frac{\psi(r)}{B} + \frac{1}{B} \sum_{j=0}^{k} \lambda^{-j} \psi(r\lambda^{j+1}).$$

By our assumption it follows that for $\beta \geq 1$, $\psi(\beta r) \leq \sqrt{\beta}\psi(r)$. Hence,

$$\mathbb{E}R_n\mathcal{G}_r \leq \frac{1}{B}\psi(r)\left(1 + \sqrt{\lambda}\sum_{j=0}^{k} \lambda^{-j/2}\right),$$

and taking $\lambda = 4$, the right-hand side is upper bounded by $5\psi(r)/B$. Moreover, for $r \geq r^*$, $\psi(r) \leq \sqrt{r/r^*}\psi(r^*) = \sqrt{rr^*}$, and thus

$$V_r^+ \leq \frac{10(1+\alpha)}{B}\sqrt{rr^*} + \sqrt{\frac{2rx}{n}} + (b-a)\left(\frac{1}{3} + \frac{1}{\alpha}\right)\frac{x}{n}.$$

Set $A = 10(1+\alpha)\sqrt{r^*}/B + \sqrt{2x/n}$ and $C = (b-a)(1/3 + 1/\alpha)x/n$, and note that $V_r^+ \leq A\sqrt{r} + C$.

We now show that $r$ can be chosen such that $V_r^+ \leq r/(\lambda BK)$. Indeed, consider the largest solution $r_0$ of $A\sqrt{r} + C = r/(\lambda BK)$. It satisfies $r_0 \geq \lambda^2 A^2 B^2 K^2/2 \geq r^*$ and $r_0 \leq (\lambda BK)^2 A^2 + 2\lambda BKC$, so that applying Lemma 3.8, it follows that every $f \in \mathcal{F}$ satisfies

$$Pf \leq \frac{K}{K-1}P_n f + \lambda BKA^2 + 2C$$

$$= \frac{K}{K-1}P_n f + \lambda BK\left(100(1+\alpha)^2 r^*/B^2 + \frac{20(1+\alpha)}{B}\sqrt{\frac{2xr^*}{n}} + \frac{2x}{n}\right)$$

$$+ (b-a)\left(\frac{1}{3} + \frac{1}{\alpha}\right)\frac{x}{n}.$$

Setting $\alpha = 1/10$ and using Lemma A.3 to show that $\sqrt{2xr^*/n} \leq Bx/(5n) + 5r^*/(2B)$ completes the proof of the first statement. The second statement is proved in the same way, by considering $V_r^-$ instead of $V_r^+$. □

PROOF OF THEOREM 3.3, SECOND PART. The proof of this result uses the same argument as for the first part. However, we consider the class $\tilde{\mathcal{G}}_r$ defined above. One can easily check that $\tilde{\mathcal{G}}_r \subset \{f \in \text{star}(\mathcal{F}, 0) : T(f) \leq r\}$, and thus $\mathbb{E}R_n\tilde{\mathcal{G}}_r \leq \psi(r)/B$. Applying Theorem 2.1 to $\tilde{\mathcal{G}}_r$, it follows that, for all $x > 0$, with probability $1 - e^{-x}$,

$$\tilde{V}_r^+ \leq \frac{2(1+\alpha)}{B}\psi(r) + \sqrt{\frac{2rx}{n}} + (b-a)\left(\frac{1}{3} + \frac{1}{\alpha}\right)\frac{x}{n}.$$

The reasoning is then the same as for the first part, and we use in the very last step that $\sqrt{2xr^*/n} \leq Bx/n + r^*/(2B)$, which gives the displayed constants. □

PROOF OF LEMMA 3.6. The map $\alpha \mapsto \alpha^2$ is Lipschitz with constant 2 when $\alpha$ is restricted to $[-1, 1]$. Applying Theorem A.6,

$$(3.2) \qquad r \geq 10\mathbb{E}R_n\{f^2 : f \in \text{star}(\mathcal{F}, 0), Pf^2 \leq r\} + \frac{26x}{n}.$$

Clearly, if $f \in \mathcal{F}$, then $f^2$ maps to $[0, 1]$ and $\text{Var}[f^2] \leq Pf^2$. Thus, Theorem 2.1 can be applied to the class $\mathcal{G}_r = \{rf^2/(Pf^2 \vee r) : f \in \mathcal{F}\}$, whose functions have range in $[0, 1]$ and variance bounded by $r$. Therefore, with probability at least $1 - e^{-x}$, every $f \in \mathcal{F}$ satisfies

$$r\frac{Pf^2 - P_nf^2}{Pf^2 \vee r} \leq 2(1+\alpha)\mathbb{E}R_n\mathcal{G}_r + \sqrt{\frac{2rx}{n}} + \left(\frac{1}{3} + \frac{1}{\alpha}\right)\frac{x}{n}.$$

Select $\alpha = 1/4$ and notice that $\sqrt{2rx/n} \leq r/4 + 2x/n$ to get

$$r\frac{Pf^2 - P_nf^2}{Pf^2 \vee r} \leq \frac{5}{2}\mathbb{E}R_n\mathcal{G}_r + \frac{r}{2} + \frac{19x}{3n}.$$

Hence, one either has $Pf^2 \leq r$, or when $Pf^2 \geq r$, since it was assumed that $P_nf^2 \leq r$,

$$Pf^2 \leq r + \frac{Pf^2}{r}\left(\frac{5}{2}\mathbb{E}R_n\mathcal{G}_r + \frac{r}{4} + \frac{19x}{3n}\right).$$

Now, if $g \in \mathcal{G}_r$, there exists $f_0 \in \mathcal{F}$ such that $g = rf_0^2/(Pf_0^2 \vee r)$. If $Pf_0^2 \leq r$, then $g = f_0^2$. On the other hand, if $Pf_0^2 > r$, then $g = rf_0^2/Pf_0^2 = f_1^2$ with $f_1 \in \text{star}(\mathcal{F}, 0)$ and $Pf_1^2 \leq r$, which shows that

$$\mathbb{E}R_n\mathcal{G}_r \leq \mathbb{E}R_n\{f^2 : f \in \text{star}(\mathcal{F}, 0), Pf^2 \leq r\}.$$

Thus, by (3.2), $Pf^2 \leq 2r$, which concludes the proof. □

**4. Data-dependent error bounds.** The results presented thus far use distribution-dependent measures of complexity of the class at hand. Indeed, the sub-root function $\psi$ of Theorem 3.3 is bounded in terms of the Rademacher averages of the star-hull of $\mathcal{F}$, but these averages can only be computed if one knows the distribution $P$. Otherwise, we have seen that it is possible to compute an upper bound on the Rademacher averages using a priori global or distribution-free knowledge about the complexity of the class at hand (such as the VC-dimension). In this section we present error bounds that can be computed directly from the data, without a priori information. Instead of computing $\psi$, we compute an estimate, $\widehat{\psi}_n$, of it. The function $\widehat{\psi}_n$ is defined using the data and is an upper bound on $\psi$ with high probability.

To simplify the exposition we restrict ourselves to the case where the functions have a range which is symmetric around zero, say $[-1, 1]$. Moreover, we can only treat the special case where $T(f) = Pf^2$, but this is a minor restriction as in most applications this is the function of interest [i.e., for which one can show $T(f) \leq BPf$].

4.1. *Results.* We now present the main result of this section, which gives an analogue of the second part of Theorem 3.3, with a completely empirical bound (i.e., the bound can be computed from the data only).

THEOREM 4.1. *Let $\mathcal{F}$ be a class of functions with ranges in $[-1, 1]$ and assume that there is some constant $B$ such that for every $f \in \mathcal{F}$, $Pf^2 \leq BPf$. Let $\widehat{\psi}_n$ be a sub-root function and let $\hat{r}^*$ be the fixed point of $\widehat{\psi}_n$. Fix $x > 0$ and assume that $\widehat{\psi}_n$ satisfies, for any $r \geq \hat{r}^*$,*

$$\widehat{\psi}_n(r) \geq c_1 \mathbb{E}_\sigma R_n\{f \in \mathrm{star}(\mathcal{F}, 0) : P_n f^2 \leq 2r\} + \frac{c_2 x}{n},$$

*where $c_1 = 2(10 \vee B)$ and $c_2 = c_1 + 11$. Then, for any $K > 1$ with probability at least $1 - 3e^{-x}$,*

$$\forall f \in \mathcal{F} \qquad Pf \leq \frac{K}{K-1} P_n f + \frac{6K}{B} \hat{r}^* + \frac{x(11 + 5BK)}{n}.$$

*Also, with probability at least $1 - 3e^{-x}$,*

$$\forall f \in \mathcal{F} \qquad P_n f \leq \frac{K+1}{K} P f + \frac{6K}{B} \hat{r}^* + \frac{x(11 + 5BK)}{n}.$$

Although these are data-dependent bounds, they are not necessarily easy to compute. There are, however, favorable interesting situations where they can be computed efficiently, as Section 6 shows.

It is natural to wonder how close the quantity $\hat{r}^*$ appearing in the above theorem is to the quantity $r^*$ of Theorem 3.3. The next theorem shows that they are close with high probability.

THEOREM 4.2. *Let $\mathcal{F}$ be a class of functions with ranges in $[-1,1]$. Fix $x > 0$ and consider the sub-root functions*

$$\psi(r) = \mathbb{E}R_n\{f \in \text{star}(\mathcal{F},0) : Pf^2 \leq r\}$$

*and*

$$\widehat{\psi}_n(r) = c_1 \mathbb{E}_\sigma R_n\{f \in \text{star}(\mathcal{F},0) : P_n f^2 \leq 2r\} + \frac{c_2 x}{n},$$

*with fixed points $r^*$ and $\hat{r}^*$, respectively, and with $c_1 = 2(10 \vee B)$ and $c_2 = 13$. Assume that $r^* \geq c_3 x/n$, where $c_3 = 26 \vee (c_2 + 2c_1)/3$. Then, with probability at least $1 - 4e^{-x}$,*

$$r^* \leq \hat{r}^* \leq 9(1+c_1)^2 r^*.$$

Thus, with high probability, $\hat{r}^*$ is an upper bound on $r^*$ and has the same asymptotic behavior. Notice that there was no attempt to optimize the constants in the above theorem. In addition, the constant $9(1+c_1)^2$ (equal to 3969 if $B \leq 10$) in Theorem 4.2 does not appear in the upper bound of Theorem 4.1.

4.2. *Proofs.* The idea of the proofs is to show that one can upper bound $\psi$ by an empirical estimate (with high probability). This requires two steps: the first one uses the concentration of the Rademacher averages to upper bound the expected Rademacher averages by their empirical versions. The second step uses Corollary 2.2 to prove that the ball over which the averages are computed [which is an $L_2(P)$ ball] can be replaced by an empirical one. Thus, $\widehat{\psi}_n$ is an upper bound on $\psi$, and one can apply Theorem 3.3, together with the following lemma, which shows how fixed points of sub-root functions relate when the functions are ordered.

LEMMA 4.3. *Suppose that $\psi, \widehat{\psi}_n$ are sub-root. Let $r^*$ (resp. $\hat{r}^*$) be the fixed point of $\psi$ (resp. $\widehat{\psi}_n$). If for $0 \leq \alpha \leq 1$ we have $\alpha\widehat{\psi}_n(r^*) \leq \psi(r^*) \leq \widehat{\psi}_n(r^*)$, then*

$$\alpha^2 \hat{r}^* \leq r^* \leq \hat{r}^*.$$

PROOF. Denoting by $\hat{r}^*_\alpha$ the fixed point of the sub-root function $\alpha\widehat{\psi}_n$, then, by Lemma 3.2 $\hat{r}^*_\alpha \leq r^* \leq \hat{r}^*$. Also, since $\widehat{\psi}_n$ is sub-root, $\widehat{\psi}_n(\alpha^2 \hat{r}^*) \geq \alpha\widehat{\psi}_n(\hat{r}^*) = \alpha\hat{r}^*$, which means $\alpha\widehat{\psi}_n(\alpha^2 \hat{r}^*) \geq \alpha^2 \hat{r}^*$. Hence, Lemma 3.2 yields $\hat{r}^*_\alpha \geq \alpha^2 \hat{r}^*$. □

PROOF OF THEOREM 4.1. Consider the sub-root function

$$\psi_1(r) = \frac{c_1}{2}\mathbb{E}R_n\{f \in \text{star}(\mathcal{F},0) : Pf^2 \leq r\} + \frac{(c_2-c_1)x}{n},$$

with fixed point $r_1^*$. Applying Corollary 2.2 when $r \geq \psi_1(r)$, it follows that with probability at least $1 - e^{-x}$,

$$\{f \in \text{star}(\mathcal{F}, 0) : Pf^2 \leq r\} \subseteq \{f \in \text{star}(\mathcal{F}, 0) : P_n f^2 \leq 2r\}.$$

Using this together with the first inequality of Lemma A.4 (with $\alpha = 1/2$) shows that if $r \geq \psi_1(r)$, with probability at least $1 - 2e^{-x}$,

$$\psi_1(r) = \frac{c_1}{2} \mathbb{E} R_n \{f \in \text{star}(\mathcal{F}, 0) : Pf^2 \leq r\} + \frac{(c_2 - c_1)x}{n}$$

$$\leq c_1 \mathbb{E}_\sigma R_n \{f \in \text{star}(\mathcal{F}, 0) : Pf^2 \leq r\} + \frac{c_2 x}{n}$$

$$\leq c_1 \mathbb{E}_\sigma R_n \{f \in \text{star}(\mathcal{F}, 0) : P_n f^2 \leq 2r\} + \frac{c_2 x}{n}$$

$$\leq \widehat{\psi}_n(r).$$

Choosing $r = r_1^*$, Lemma 4.3 shows that with probability at least $1 - 2e^{-x}$,

(4.1) $$r_1^* \leq \hat{r}^*.$$

Also, for all $r \geq 0$,

$$\psi_1(r) \geq B \mathbb{E} R_n \{f \in \text{star}(\mathcal{F}, 0) : Pf^2 \leq r\},$$

and so from Theorem 3.3, with probability at least $1 - e^{-x}$, every $f \in \mathcal{F}$ satisfies

$$Pf \leq \frac{K}{K-1} P_n f + \frac{6K r_1^*}{B} + \frac{(11 + 5BK)x}{n}.$$

Combining this with (4.1) gives the first result. The second result is proved in a similar manner. $\square$

PROOF OF THEOREM 4.2. Consider the functions

$$\psi_1(r) = \frac{c_1}{2} \mathbb{E} R_n \{f \in \text{star}(\mathcal{F}, 0) : Pf^2 \leq r\} + \frac{(c_2 - c_1)x}{n}$$

and

$$\psi_2(r) = c_1 \mathbb{E} R_n \{f \in \text{star}(\mathcal{F}, 0) : Pf^2 \leq r\} + \frac{c_3 x}{n},$$

and denote by $r_1^*$ and $r_2^*$ the fixed points of $\psi_1$ and $\psi_2$, respectively. The proof of Theorem 4.1 shows that with probability at least $1 - 2e^{-x}$, $r_1^* \leq \hat{r}^*$.

Now apply Lemma 3.6 to show that if $r \geq \psi_2(r)$, then with probability at least $1 - e^{-x}$,

$$\{f \in \text{star}(\mathcal{F}, 0) : P_n f^2 \leq r\} \subseteq \{f \in \text{star}(\mathcal{F}, 0) : Pf^2 \leq 2r\}.$$

Using this together with the second inequality of Lemma A.4 (with $\alpha = 1/2$) shows that if $r \geq \psi_2(r)$, with probability at least $1 - 2e^{-x}$,

$$\begin{aligned}
\widehat{\psi}_n(r) &= c_1 \mathbb{E}_\sigma R_n\{f \in \text{star}(\mathcal{F}, 0) : P_n f^2 \leq 2r\} + \frac{c_2 x}{n} \\
&\leq c_1 \sqrt{2} \mathbb{E}_\sigma R_n\{f \in \text{star}(\mathcal{F}, 0) : P_n f^2 \leq r\} + \frac{c_2 x}{n} \\
&\leq c_1 \sqrt{2} \mathbb{E}_\sigma R_n\{f \in \text{star}(\mathcal{F}, 0) : P f^2 \leq 2r\} + \frac{c_2 x}{n} \\
&\leq \frac{3\sqrt{2}}{2} c_1 \mathbb{E} R_n\{f \in \text{star}(\mathcal{F}, 0) : P f^2 \leq 2r\} + \frac{(c_2 + 2c_1)x}{n} \\
&\leq 3 c_1 \mathbb{E} R_n\{f \in \text{star}(\mathcal{F}, 0) : P f^2 \leq r\} + \frac{(c_2 + 2c_1)x}{n} \\
&\leq 3 \psi_2(r),
\end{aligned}$$

where the sub-root property was used twice (in the first and second to last inequalities). Lemma 4.3 thus gives $\hat{r}^* \leq 9 r_2^*$.

Also notice that for all $r$, $\psi(r) \leq \psi_1(r)$, and hence $r^* \leq r_1^*$. Moreover, for all $r \geq \psi(r)$ (hence $r \geq r^* \geq c_3 x/n$), $\psi_2(r) \leq c_1 \psi(r) + r$, so that $\psi_2(r^*) \leq (c_1 + 1) r^* = (c_1 + 1) \psi(r^*)$. Lemma 4.3 implies that $r_2^* \leq (1 + c_1)^2 r^*$. □

**5. Prediction with bounded loss.** In this section we discuss the application of our results to prediction problems, such as classification and regression. For such problems there are an *input space* $\mathcal{X}$ and an *output space* $\mathcal{Y}$, and the product $\mathcal{X} \times \mathcal{Y}$ is endowed with an unknown probability measure $P$. For example, classification corresponds to the case where $\mathcal{Y}$ is discrete, typically $\mathcal{Y} = \{-1, 1\}$, and regression corresponds to the continuous case, typically $\mathcal{Y} = [-1, 1]$. Note that assuming the boundedness of the target values is a typical assumption in theoretical analysis of regression procedures. To analyze the case of unbounded targets, one usually truncates the values at a certain threshold and bounds the probability of exceeding that threshold (see, e.g., the techniques developed in [12]).

The *training sample* is a sequence $(X_1, Y_1), \ldots, (X_n, Y_n)$ of $n$ independent and identically distributed (i.i.d.) pairs sampled according to $P$. A loss function $\ell : \mathcal{Y} \times \mathcal{Y} \to [0, 1]$ is defined and the goal is to find a function $f : \mathcal{X} \to \mathcal{Y}$ from a class $\mathcal{F}$ that minimizes the expected loss

$$\mathbb{E} \ell_f = \mathbb{E} \ell(f(X), Y).$$

Since the probability distribution $P$ is unknown, one cannot directly minimize the expected loss over $\mathcal{F}$.

The key property that is needed to apply our results is the fact that $\text{Var}[f] \leq BPf$ (or $Pf^2 \leq BPf$ to obtain data-dependent bounds). This will trivially be the case for the class $\{\ell_f : f \in \mathcal{F}\}$, as all its functions are uniformly bounded and nonnegative. This case, studied in Section 5.1, is, however, not the most interesting. Indeed, it is when one studies the excess risk $\ell_f - \ell_{f^*}$ that our approach shows its superiority over previous ones; when the class $\{\ell_f - \ell_{f^*}\}$ satisfies the variance condition (and Section 5.2 gives examples of this), we obtain distribution-dependent bounds that are optimal in certain cases, and data-dependent bounds of the same order.

5.1. *General results without assumptions.* Define the following class of functions, called the *loss class associated with* $\mathcal{F}$:

$$\ell_\mathcal{F} = \{\ell_f : f \in \mathcal{F}\} = \{(x,y) \mapsto \ell(f(x), y) : f \in \mathcal{F}\}.$$

Notice that $\ell_\mathcal{F}$ is a class of nonnegative functions. Applying Theorem 4.1 to this class of functions gives the following corollary.

COROLLARY 5.1. *For a loss function* $\ell : \mathcal{Y} \times \mathcal{Y} \to [0,1]$, *define*

$$\widehat{\psi}_n(r) = 20 \mathbb{E}_\sigma R_n \{f \in \text{star}(\ell_\mathcal{F}, 0) : P_n f^2 \leq 2r\} + \frac{13x}{n},$$

*with fixed point* $\hat{r}^*$. *Then, for any* $K > 1$ *with probability at least* $1 - 3e^{-x}$,

$$\forall f \in \mathcal{F} \qquad P\ell_f \leq \frac{K}{K-1} P_n \ell_f + 6K \hat{r}^* + \frac{x(11+5K)}{n}.$$

A natural approach is to minimize the empirical loss $P_n \ell_f$ over the class $\mathcal{F}$. The following result shows that this approach leads to an estimate with expected loss near minimal. How close it is to the minimal expected loss depends on the value of the minimum, as well as on the local Rademacher averages of the class.

THEOREM 5.2. *For a loss function* $\ell : \mathcal{Y} \times \mathcal{Y} \to [0,1]$, *define* $\psi(r)$, $\widehat{\psi}_n(r)$, $r^*$ *and* $\hat{r}^*$ *as in Theorem* 5.1. *Let* $L^* = \inf_{f \in \mathcal{F}} P\ell_f$. *Then there is a constant* $c$ *such that with probability at least* $1 - 2e^{-x}$, *the minimizer* $\hat{f} \in \mathcal{F}$ *of* $P_n \ell_f$ *satisfies*

$$P\ell_{\hat{f}} \leq L^* + c(\sqrt{L^* r^*} + r^*).$$

*Also, with probability at least* $1 - 4e^{-x}$,

$$P\ell_{\hat{f}} \leq L^* + c(\sqrt{L^* \hat{r}^*} + \hat{r}^*).$$

The proof of this theorem is given in Appendix A.2.

This theorem has the same flavor as Theorem 4.2 of [19]. We have not used any property besides the positivity of the functions in the class. This indicates that there might not be a significant gain compared to earlier results (as without further assumptions the optimal rates are known). Indeed, a careful examination of this result shows that when $L^* > 0$, the difference between $P\ell_{\hat{f}}$ and $L^*$ is essentially of order $\sqrt{r^*}$. For a class of $\{0,1\}$-valued functions with VC-dimension $d$, for example, this would be $\sqrt{d \log n / n}$. On the other hand, the result of [19] is more refined since the Rademacher averages are not localized around 0 (as they are here), but rather around the minimizer of the empirical error itself. Unfortunately, the small ball in [19] is not defined as $P_n \ell_f \leq P_n \ell_{\hat{f}} + r$ but as $P_n \ell_f \leq 16 P_n \ell_{\hat{f}} + r$. This means that in the general situation where $L^* > 0$, since $P_n \ell_{\hat{f}}$ does not converge to 0 with increasing $n$ (as it is expected to be close to $P\ell_{\hat{f}}$ which itself converges to $L^*$), the radius of the ball around $\ell_{\hat{f}}$ (which is $15 P_n \ell_{\hat{f}} + r$) will not converge to 0. Thus, the localized Rademacher average over this ball will converge at speed $\sqrt{d/n}$. In other words, our Theorem 5.2 and Theorem 4.2 of [19] essentially have the same behavior. But this is not surprising, as it is known that this is the optimal rate of convergence in this case. To get an improvement in the rates of convergence, one needs to make further assumptions on the distribution $P$ or on the class $\mathcal{F}$.

5.2. *Improved results for the excess risk.* Consider a loss function $\ell$ and function class $\mathcal{F}$ that satisfy the following conditions.

1. For every probability distribution $P$ there is an $f^* \in \mathcal{F}$ satisfying $P\ell_{f^*} = \inf_{f \in \mathcal{F}} P\ell_f$.
2. There is a constant $L$ such that $\ell$ is $L$-Lipschitz in its first argument: for all $y, \hat{y}_1, \hat{y}_2$,
$$|\ell(\hat{y}_1, y) - \ell(\hat{y}_2, y)| \leq L|\hat{y}_1 - \hat{y}_2|.$$
3. There is a constant $B \geq 1$ such that for every probability distribution and every $f \in \mathcal{F}$,
$$P(f - f^*)^2 \leq B P(\ell_f - \ell_{f^*}).$$

These conditions are not too restrictive as they are met by several commonly used regularized algorithms with convex losses.

Note that condition 1 could be weakened, and one could consider a function which is only close to achieving the infimum, with an appropriate change to condition 3. This generalization is straightforward, but it would make the results less readable, so we omit it.

Condition 2 implies that, for all $f \in \mathcal{F}$,
$$P(\ell_f - \ell_{f^*})^2 \leq L^2 P(f - f^*)^2.$$

Condition 3 usually follows from a uniform convexity condition on $\ell$. An important example is the quadratic loss $\ell(y, y') = (y - y')^2$, when the function class $\mathcal{F}$ is convex and uniformly bounded. In particular, if $|f(x) - y| \in [0, 1]$ for all $f \in \mathcal{F}$, $x \in \mathcal{X}$ and $y \in \mathcal{Y}$, then the conditions are satisfied with $L = 2$ and $B = 1$ (see [18]). Other examples are described in [26] and in [2].

The first result we present is a direct but instructive corollary of Theorem 3.3.

COROLLARY 5.3. *Let $\mathcal{F}$ be a class of functions with ranges in $[-1, 1]$ and let $\ell$ be a loss function satisfying conditions 1–3 above. Let $\hat{f}$ be any element of $\mathcal{F}$ satisfying $P_n \ell_{\hat{f}} = \inf_{f \in \mathcal{F}} P_n \ell_f$. Assume $\psi$ is a sub-root function for which*

$$\psi(r) \geq BL\mathbb{E}R_n\{f \in \mathcal{F} : L^2 P(f - f^*)^2 \leq r\}.$$

*Then for any $x > 0$ and any $r \geq \psi(r)$, with probability at least $1 - e^{-x}$,*

$$P(\ell_{\hat{f}} - \ell_{f^*}) \leq 705 \frac{r}{B} + \frac{(11L + 27B)x}{n}.$$

PROOF. One applies Theorem 3.3 (first part) to the class $\ell_f - \ell_{f^*}$ with $T(f) = L^2 P(f - f^*)^2$ and uses the fact that by Theorem A.6, and by the symmetry of the Rademacher variables, $L\mathbb{E}R_n\{f : L^2 P(f - f^*)^2 \leq r\} \geq \mathbb{E}R_n\{\ell_f - \ell_{f^*} : L^2 P(f - f^*)^2 \leq r\}$. The result follows from noticing that $P_n(\ell_{\hat{f}} - \ell_{f^*}) \leq 0$. □

Instead of comparing the loss of $f$ to that of $f^*$, one could compare it to the loss of the best measurable function (the regression function for regression function estimation, or the Bayes classifier for classification). The techniques proposed here can be adapted to this case.

Using Corollary 5.3, one can (with minor modification) recover the results of [22] for model selection. These have been shown to match the minimax results in various situations. In that sense, Corollary 5.3 can be considered as sharp.

Next we turn to the main result of this section. It is a version of Corollary 5.3 with a fully data-dependent bound. This is obtained by modifying $\psi$ in three ways: the Rademacher averages are replaced by empirical ones, the radius of the ball is in the $L_2(P_n)$ norm instead of $L_2(P)$, and finally, the center of the ball is $\hat{f}$ instead of $f^*$.

THEOREM 5.4. *Let $\mathcal{F}$ be a convex class of functions with range in $[-1, 1]$ and let $\ell$ be a loss function satisfying conditions 1–3 above. Let $\hat{f}$ be any element of $\mathcal{F}$ satisfying $P_n \ell_{\hat{f}} = \inf_{f \in \mathcal{F}} P_n \ell_f$. Define*

$$(5.1) \qquad \widehat{\psi}_n(r) = c_1 \mathbb{E}_\sigma R_n\{f \in \mathcal{F} : P_n(f - \hat{f})^2 \leq c_3 r\} + \frac{c_2 x}{n},$$

where $c_1 = 2L(B \vee 10L)$, $c_2 = 11L^2 + c_1$ and $c_3 = 2824 + 4B(11L + 27B)/c_2$. Then with probability at least $1 - 4e^{-x}$,

$$P(\ell_{\hat{f}} - \ell_{f^*}) \leq \frac{705}{B}\hat{r}^* + \frac{(11L + 27B)x}{n},$$

where $\hat{r}^*$ is the fixed point of $\widehat{\psi}_n$.

REMARK 5.5. Unlike Corollary 5.3, the class $\mathcal{F}$ in Theorem 5.4 has to be convex. This ensures that it is star-shaped around any of its elements (which implies that $\widehat{\psi}_n$ is sub-root even though $\hat{f}$ is random). However, convexity of the loss class is not necessary, so that this theorem still applies to many situations of interest, in particular to regularized regression, where the functions are taken in a vector space or a ball of a vector space.

REMARK 5.6. Although the theorem is stated with explicit constants, there is no reason to think that these are optimal. The fact that the constant 705 appears actually is due to our failure to apply the second part of Theorem 3.3 to the initial loss class, which is not star-shaped (this would have given a 7 instead). However, with some additional effort, one can probably obtain much better constants.

As we explained earlier, although the statement of Theorem 5.4 is similar to Theorem 4.2 in [19], there is an important difference in the way the localized averages are defined: in our case the radius is a constant times $r$, while in [19] there is an additional term, involving the loss of the empirical risk minimizer, which may not converge to zero. Hence, the complexity decreases faster in our bound.

The additional property required in the proof of this result compared to the proof of Theorem 4.1 is that under the assumptions of the theorem, the minimizers of the empirical loss and of the true loss are close with respect to the $L_2(P)$ and the $L_2(P_n)$ distances (this has also been used in [20] and [31, 32]).

PROOF OF THEOREM 5.4. Define the function $\psi$ as

(5.2) $\qquad \psi(r) = \frac{c_1}{2}\mathbb{E}R_n\{f \in \mathcal{F} : L^2 P(f - f^*)^2 \leq r\} + \frac{(c_2 - c_1)x}{n}.$

Notice that since $\mathcal{F}$ is convex and thus star-shaped around each of its points, Lemma 3.4 implies that $\psi$ is sub-root. Now, for $r \geq \psi(r)$ Corollary 5.3 and condition 3 on the loss function imply that, with probability at least $1 - e^{-x}$,

(5.3) $\quad L^2 P(\hat{f} - f^*)^2 \leq BL^2 P(\ell_{\hat{f}} - \ell_{f^*}) \leq 705L^2 r + \frac{(11L + 27B)BL^2 x}{n}.$

Denote the right-hand side by $s$. Since $s \geq r \geq r^*$, then $s \geq \psi(s)$ (by Lemma 3.2), and thus

$$s \geq 10L^2 \mathbb{E} R_n \{f \in \mathcal{F} : L^2 P(f - f^*)^2 \leq s\} + \frac{11L^2 x}{n}.$$

Therefore, Corollary 2.2 applied to the class $L\mathcal{F}$ yields that with probability at least $1 - e^{-x}$,

$$\{f \in \mathcal{F}, L^2 P(f - f^*)^2 \leq s\} \subset \{f \in \mathcal{F}, L^2 P_n(f - f^*)^2 \leq 2s\}.$$

This, combined with (5.3), implies that with probability at least $1 - 2e^{-x}$,

(5.4)
$$P_n(\hat{f} - f^*)^2 \leq 2\left(705r + \frac{(11L + 27B)Bx}{n}\right)$$
$$\leq 2\left(705 + \frac{(11L + 27B)B}{c_2}\right)r,$$

where the second inequality follows from $r \geq \psi(r) \geq c_2 x/n$. Define $c = 2(705 + (11L + 27B)B/c_2)$. By the triangle inequality in $L_2(P_n)$, if (5.4) occurs, then any $f \in \mathcal{F}$ satisfies

$$P_n(f - \hat{f})^2 \leq (\sqrt{P_n(f - f^*)^2} + \sqrt{P_n(f^* - \hat{f})^2})^2$$
$$\leq (\sqrt{P_n(f - f^*)^2} + \sqrt{cr})^2.$$

Appealing again to Corollary 2.2 applied to $L\mathcal{F}$ as before, but now for $r \geq \psi(r)$, it follows that with probability at least $1 - 3e^{-x}$,

$$\{f \in \mathcal{F} : L^2 P(f - f^*)^2 \leq r\}$$
$$\subseteq \{f \in \mathcal{F} : L^2 P_n(f - \hat{f})^2 \leq (\sqrt{2} + \sqrt{c})^2 L^2 r\}.$$

Combining this with Lemma A.4 shows that, with probability at least $1 - 4e^{-x}$,

$$\psi(r) \leq c_1 \mathbb{E}_\sigma R_n \{f \in \mathcal{F} : L^2 P(f - f^*)^2 \leq r\} + \frac{c_2 x}{n}$$
$$\leq c_1 \mathbb{E}_\sigma R_n \{f : P_n(f - f^*)^2 \leq (\sqrt{2} + \sqrt{c})^2 r\} + \frac{c_2 x}{n}$$
$$\leq c_1 \mathbb{E}_\sigma R_n \{f : P_n(f - f^*)^2 \leq (4 + 2c)r\} + \frac{c_2 x}{n}$$
$$\leq \widehat{\psi}_n(r).$$

Setting $r = r^*$ in the above argument and applying Lemma 4.3 shows that $r^* \leq \hat{r}^*$, which, together with (5.3), concludes the proof. $\square$

**6. Computing local Rademacher complexities.** In this section we deal with the computation of local Rademacher complexities and their fixed points. We first propose a simple iterative procedure for estimating the fixed point of an arbitrary sub-root function and then give two examples of situations where it is possible to compute an upper bound on the local Rademacher complexities. In the case of classification with the discrete loss, this can be done by solving a weighted error minimization problem. In the case of kernel classes, it is obtained by computing the eigenvalues of the empirical Gram matrix.

6.1. *The iterative procedure.* Recall that Theorem 4.1 indicates that one can obtain an upper bound in terms of empirical quantities only. However, it remains to be explained how to compute these quantities effectively. We propose to use a procedure similar to that of Koltchinskii and Panchenko [16], by applying the sub-root function iteratively. The next lemma shows that applying the sub-root function iteratively gives a sequence that converges monotonically and quickly to the fixed point.

LEMMA 6.1. *Let $\psi:[0,\infty) \to [0,\infty)$ be a (nontrivial) sub-root function. Fix $r_0 \geq r^*$, and for all $k > 0$ define $r_{k+1} = \psi(r_k)$. Then for all $N > 0$, $r_{N+1} \leq r_N$, and*

$$r^* \leq r_N \leq \left(\frac{r_0}{r^*}\right)^{2^{-N}} r^*.$$

*In particular, for any $\varepsilon > 0$, if $N$ satisfies*

$$N \geq \log_2\left(\frac{\ln(r_0/r^*)}{\ln(1+\varepsilon)}\right),$$

*then $r_N \leq (1+\varepsilon)r^*$.*

PROOF. Notice that if $r_k \geq r^*$, then $r_{k+1} = \psi(r_k) \geq \psi(r^*) = r^*$. Also,

$$\frac{\psi(r_k)}{\sqrt{r_k}} \leq \frac{\psi(r^*)}{\sqrt{r^*}} = \sqrt{r^*} \leq \sqrt{r_k},$$

and so $r_{k+1} \leq r_k$ and $r_{k+1}/r^* \leq (r_k/r^*)^{1/2}$. An easy induction shows that $r_N/r^* \leq (r_0/r^*)^{2^{-N}}$. □

Notice that in the results of [16], the analysis of the iterative procedure was tied to the probabilistic upper bounds. However, here we make the issues separate: the bounds of previous sections are valid no matter how the fixed point is estimated. In the above lemma, one can use a random sub-root function.

6.2. *Local Rademacher complexities for classification loss classes.* Consider the case where $\mathcal{Y} = \{-1, 1\}$ and the loss is the discrete loss, $\ell(y, y') = \mathbf{1}[y \neq y']$. Since $\ell^2 = \ell$, one can write

$$\mathbb{E}_\sigma R_n\{f \in \text{star}(\ell_\mathcal{F}, 0) : P_n f^2 \leq 2r\}$$
$$= \mathbb{E}_\sigma R_n\{\alpha \ell_f : \alpha \in (0, 1], f \in \mathcal{F}, P_n \ell_f^2 \leq 2r/\alpha^2\}$$
$$= \mathbb{E}_\sigma R_n\{\alpha \ell_f : \alpha \in (0, 1], f \in \mathcal{F}, P_n \ell_f \leq 2r/\alpha^2\}$$
$$= \sup_{\alpha \in (0,1]} \alpha \mathbb{E}_\sigma R_n\{\ell_f : f \in \mathcal{F}, P_n \ell_f \leq 2r/\alpha^2\}$$
$$= \sup_{\alpha \in [\sqrt{2r},1]} \alpha \mathbb{E}_\sigma R_n\{\ell_f : f \in \mathcal{F}, P_n \ell_f \leq 2r/\alpha^2\},$$

where the last equality follows from the fact that $P_n \ell_f \leq 1$ for all $f$. Substituting into Corollary 5.1 gives the following result.

COROLLARY 6.2. *Let $\mathcal{Y} = \{\pm 1\}$, let $\ell$ be the discrete loss defined on $\mathcal{Y}$ and let $\mathcal{F}$ be a class of functions with ranges in $\mathcal{Y}$. Fix $x > 0$ and define*

$$\widehat{\psi}_n(r) = 20 \sup_{\alpha \in [\sqrt{2r},1]} \alpha \mathbb{E}_\sigma R_n\{\ell_f : f \in \mathcal{F}, P_n \ell_f \leq 2r/\alpha^2\} + \frac{26x}{n}.$$

*Then for all $K > 1$, with probability at least $1 - 3e^{-x}$, for all $f \in \mathcal{F}$,*

$$P\ell_f \leq \frac{K}{K-1} P_n \ell_f + cK\left(\hat{r}^* + \frac{x}{n}\right),$$

*where $\hat{r}^*$ is the fixed point of $\widehat{\psi}_n$.*

The following theorem shows that upper bounds on $\widehat{\psi}_n(r)$ can by computed whenever one can perform weighted empirical risk minimization. In other words, if there is an efficient algorithm for minimizing a weighted sum of classification errors, there is an efficient algorithm for computing an upper bound on the localized Rademacher averages. The empirical minimization algorithm needs to be run repeatedly on different realizations of the $\sigma_i$, but with fast convergence toward the expectation as the number of iterations grows. A similar result was known for global Rademacher averages and this shows that the localization and the use of star-hulls do not greatly affect the computational complexity.

THEOREM 6.3. *The empirical local Rademacher complexity of the classification loss class, defined in Corollary 6.2, satisfies*

$$\widehat{\psi}_n(r) = c \sup_{\alpha \in [\sqrt{2r},1]} \alpha \mathbb{E}_\sigma R_n\{\ell_f : f \in \mathcal{F}, P_n \ell_f \leq 2r/\alpha^2\} + \frac{26x}{n}$$

$$\leq c \sup_{\alpha \in [\sqrt{2r},1]} \alpha \mathbb{E}_\sigma \min_{\mu \geq 0} \left( \left( \frac{2r}{\alpha^2} - \frac{1}{2} \right) \mu + \frac{1}{2n} \sum_{i=1}^n |\sigma_i + \mu Y_i| - J(\mu) \right) + \frac{26x}{n},$$

*where*

$$J(\mu) = \min_{f \in \mathcal{F}} \frac{1}{n} \sum_{i=1}^n |\sigma_i + \mu Y_i| \ell(f(X_i), \text{sign}(\sigma_i + \mu Y_i)).$$

The quantity $J(\mu)$ can be viewed as the minimum of a certain weighted empirical risk when the labels are corrupted by noise and the noise level is determined by the parameter (Lagrange multiplier) $\mu$. Using the fact that $J(\mu)$ is Lipschitz in $\mu$, a finite grid of values of $J(\mu)$ can be used to obtain a function $\phi$ that is an upper bound on $\widehat{\psi}_n$. Then the function $r \mapsto \sqrt{r} \sup_{r'} \phi(r')/\sqrt{r'}$ is a sub-root upper bound on $\widehat{\psi}_n$.

In order to prove Theorem 6.3 we need the following lemma (adapted from [1]) which relates the localized Rademacher averages to a weighted error minimization problem.

LEMMA 6.4. *For every $b \in [0,1]$,*

$$\mathbb{E}_\sigma R_n\{\ell_f : f \in \mathcal{F}, P_n \ell_f \leq b\}$$
$$= \tfrac{1}{2} - \mathbb{E}_\sigma \min\{P_n \ell(f(X), \sigma) : f \in \mathcal{F}, P_n \ell(f(X), Y) \leq b\}.$$

PROOF. Notice that for $y, y' \in \{\pm 1\}$, $\ell(y, y') = \mathbf{1}[y \neq y'] = |y - y'|/2$. Thus

$$2 \sum_{i=1}^n \sigma_i \ell(f(X_i), Y_i) = \sum_{i:Y_i=1} \sigma_i |f(X_i) - 1| + \sum_{i:Y_i=-1} \sigma_i |f(X_i) + 1|$$
$$= \sum_{i:Y_i=1} \sigma_i(2 - |f(X_i) + 1|) + \sum_{i:Y_i=-1} \sigma_i |f(X_i) + 1|$$
$$= \sum_{i=1}^n -Y_i \sigma_i |f(X_i) + 1| + 2 \sum_{i:Y_i=1} \sigma_i.$$

Because of the symmetry of $\sigma_i$, for fixed $X_i$ the vector $(-Y_i \sigma_i)_{i=1}^n$ has the same distribution as $(\sigma_i)_{i=1}^n$. Thus when we take the expectation, we can replace $-Y_i \sigma_i$ by $\sigma_i$. Moreover, we have

$$\sum_{i=1}^n \sigma_i |f(X_i) + 1| = \sum_{i:\sigma_i=1} |f(X_i) + 1| + \sum_{i:\sigma_i=-1} -|f(X_i) + 1|$$
$$= \sum_{i:\sigma_i=1} (2 - |f(X_i) - 1|) + \sum_{i:\sigma_i=-1} -|f(X_i) + 1|$$
$$= \sum_{i=1}^n -|f(X_i) - \sigma_i| + 2 \sum_{i:\sigma_i=-1} 1,$$

implying that

$$\mathbb{E}_\sigma R_n\{\ell_f : f \in \mathcal{F}, P_n \ell_f \leq b\}$$

$$= \frac{1}{n}\left(\mathbb{E}_\sigma \sum_{i:Y_i=1} \sigma_i + \mathbb{E}_\sigma \sum_{i:\sigma_i=-1} 1 \right.$$

$$\left. + \mathbb{E}_\sigma \sup\{-P_n \ell(f(X), \sigma) : f \in \mathcal{F}, P_n \ell(f(X), Y) \leq b\}\right),$$

which proves the claim. □

PROOF OF THEOREM 6.3. From Lemma 6.4,

$$\widehat{\psi}_n(r) = c \sup_{\alpha \in [\sqrt{2r}, 1]} \alpha\left(\frac{1}{2} - \mathbb{E}_\sigma \min\left\{P_n \ell(f(X), \sigma) : \right.\right.$$

$$\left.\left. f \in \mathcal{F}, P_n \ell(f(X), Y) \leq \frac{2r}{\alpha^2}\right\}\right) + \frac{26x}{n}.$$

Fix a realization of the $\sigma_i$. It is easy to see that when $\mu \geq 0$, each $f$ for which $P_n \ell(f(X), Y) \leq 2r/\alpha^2$ satisfies

$$P_n \ell(f(X), \sigma) \geq P_n \ell(f(X), \sigma) + \mu\left(P_n \ell(f(X), Y) - \frac{2r}{\alpha^2}\right).$$

Let $L(f, \mu)$ denote the right-hand side and let $g(\mu) = \min_{f \in \mathcal{F}} L(f, \mu)$. Then

$$\min\{P_n \ell(f(X), \sigma) : f \in \mathcal{F}, P_n \ell(f(X), Y) \leq 2r/\alpha^2\} \geq g(\mu).$$

But, using the fact that $\ell(y, \hat{y}) = (1 - y\hat{y})/2$,

$$g(\mu) = \min_{f \in \mathcal{F}} \frac{1}{n} \sum_{i=1}^n (\ell(f(X_i), \sigma_i) + \mu \ell(f(X_i), Y_i)) - \frac{2r}{\alpha^2}$$

$$= \min_{f \in \mathcal{F}} \frac{1}{n} \sum_{i=1}^n \left(\frac{1 - f(X_i)\sigma_i}{2} + \mu \frac{1 - f(X_i)Y_i}{2}\right) - \frac{2r}{\alpha^2}$$

$$= \min_{f \in \mathcal{F}} \frac{1}{n} \sum_{i=1}^n \left(|\sigma_i + \mu Y_i| \frac{1 - f(X_i)\mathrm{sign}(\sigma_i + \mu Y_i)}{2} - \frac{|\sigma_i + \mu Y_i|}{2}\right)$$

$$+ \frac{1+\mu}{2} - \frac{2r}{\alpha^2}$$

$$= \min_{f \in \mathcal{F}} \frac{1}{n} \sum_{i=1}^n |\sigma_i + \mu Y_i| \ell(f(X_i), \mathrm{sign}(\sigma_i + \mu Y_i))$$

$$- \frac{1}{2n} \sum_{i=1}^n |\sigma_i + \mu Y_i| + \frac{1+\mu}{2} - \frac{2r}{\alpha^2}.$$

Substituting gives the result. □

LOCAL RADEMACHER COMPLEXITIES 31

6.3. *Local Rademacher complexities for kernel classes.* One case in which the functions $\psi$ and $\widehat{\psi}_n$ can be computed explicitly is when $\mathcal{F}$ is a kernel class, that is, the unit ball in the reproducing kernel Hilbert space associated with a positive definite kernel $k$. Observe that in this case $\mathcal{F}$ is a convex and symmetric set.

Let $k$ be a positive definite function on $\mathcal{X}$, that is, a symmetric function such that for all $n \geq 1$,

$$\forall x_1, \ldots, x_n \in \mathcal{X}, \ \forall \alpha_1, \ldots, \alpha_n \in \mathbb{R} \qquad \sum_{i,j=1}^{n} \alpha_i \alpha_j k(x_i, x_j) \geq 0.$$

Recall the main properties of reproducing kernel Hilbert spaces that we require:

(a) The reproducing kernel Hilbert space associated with $k$ is the unique Hilbert space $\mathcal{H}$ of functions on $\mathcal{X}$ such that for all $f \in \mathcal{F}$ and all $x \in \mathcal{X}$, $k(x, \cdot) \in \mathcal{H}$ and

(6.1) $$f(x) = \langle f, k(x, \cdot) \rangle.$$

(b) $\mathcal{H}$ can be constructed as the completion of the linear span of the functions $k(x, \cdot)$ for $x \in \mathcal{X}$, endowed with the inner product

$$\left\langle \sum_{i=1}^{n} \alpha_i k(x_i, \cdot), \sum_{j=1}^{m} \beta_j k(y_j, \cdot) \right\rangle = \sum_{i,j=1}^{n,m} \alpha_i \beta_j k(x_i, y_j).$$

We use $\|\cdot\|$ to denote the norm in $\mathcal{H}$.

One method for regression consists of solving the following least squares problem in the unit ball of $\mathcal{H}$:

$$\min_{f \in \mathcal{H} \,:\, \|f\| \leq 1} \frac{1}{n} \sum_{i=1}^{n} (f(X_i) - Y_i)^2.$$

Notice that considering a ball of some other radius is equivalent to rescaling the class. We are thus interested in computing the localized Rademacher averages of the class of functions

$$\mathcal{F} = \{f \in \mathcal{H} : \|f\| \leq 1\}.$$

Assume that $\mathbb{E}k(X, X) < \infty$ and define $T : L_2(P) \to L_2(P)$ as the integral operator associated with $k$ and $P$, that is, $Tf(\cdot) = \int k(\cdot, y) f(y) \, dP(y)$. It is possible to show that $T$ is a positive semidefinite trace-class operator. Let $(\lambda_i)_{i=1}^{\infty}$ be its eigenvalues, arranged in a nonincreasing order. Also, given an i.i.d. sample $X_1, \ldots, X_n$ from $P$, consider the normalized Gram matrix (or *kernel matrix*) $\hat{T}_n$ defined as $\hat{T}_n = \frac{1}{n}(k(X_i, X_j))_{i,j=1,\ldots,n}$. Let $(\hat{\lambda}_i)_{i=1}^{n}$ be its eigenvalues, arranged in a nonincreasing order.

The following result was proved in [24].

THEOREM 6.5. *For every $r > 0$,*

$$\mathbb{E}R_n\{f \in \mathcal{F} : Pf^2 \leq r\} \leq \left(\frac{2}{n}\sum_{i=1}^{\infty}\min\{r, \lambda_i\}\right)^{1/2}.$$

*Moreover, there exists an absolute constant $c$ such that if $\lambda_1 \geq 1/n$, then for every $r \geq 1/n$,*

$$\mathbb{E}R_n\{f \in \mathcal{F} : Pf^2 \leq r\} \geq c\left(\frac{1}{n}\sum_{i=1}^{\infty}\min\{r, \lambda_i\}\right)^{1/2}.$$

The following lemma is a data-dependent version.

LEMMA 6.6. *For every $r > 0$,*

$$\mathbb{E}_\sigma R_n\{f \in \mathcal{F} : P_n f^2 \leq r\} \leq \left(\frac{2}{n}\sum_{i=1}^{n}\min\{r, \hat{\lambda}_i\}\right)^{1/2}.$$

The proof of this result can be found in Appendix A.2. The fact that we have replaced $Pf^2$ by $P_n f^2$ and conditioned on the data yields a result that involves only the eigenvalues of the empirical Gram matrix.

We can now state a consequence of Theorem 5.4 for the proposed regression algorithm on the unit ball of $\mathcal{H}$.

COROLLARY 6.7. *Assume that $\sup_{x \in \mathcal{X}} k(x, x) \leq 1$. Let $\mathcal{F} = \{f \in \mathcal{H} : \|f\| \leq 1\}$ and let $\ell$ be a loss function satisfying conditions 1–3. Let $\hat{f}$ be any element of $\mathcal{F}$ satisfying $P_n \ell_{\hat{f}} = \inf_{f \in \mathcal{F}} P_n \ell_f$.*

*There exists a constant $c$ depending only on $L$ and $B$ such that with probability at least $1 - 6e^{-x}$,*

$$P(\ell_{\hat{f}} - \ell_{f^*}) \leq c\left(\hat{r}^* + \frac{x}{n}\right),$$

*where*

$$\hat{r}^* \leq \min_{0 \leq h \leq n}\left(\frac{h}{n} + \sqrt{\frac{1}{n}\sum_{i > h}\hat{\lambda}_i}\right).$$

We observe that $\hat{r}^*$ is at most of order $1/\sqrt{n}$ (if we take $h = 0$), but can be of order $\log n / n$ if the eigenvalues of $\hat{T}_n$ decay exponentially quickly.

In addition, the eigenvalues of the Gram matrix are not hard to compute, so that the above result can suggest an implementable heuristic for choosing the kernel $k$ from the data. The issue of the choice of the kernel is being intensively studied in the machine learning community.

PROOF. Because of the symmetry of the $\sigma_i$ and because $\mathcal{F}$ is convex and symmetric,

$$\begin{aligned}\mathbb{E}_\sigma R_n\{f \in \mathcal{F} : P_n(f - \hat{f})^2 \leq c_3 r\} &= \mathbb{E}_\sigma R_n\{f - \hat{f} : f \in \mathcal{F}, P_n(f - \hat{f})^2 \leq c_3 r\} \\ &\leq \mathbb{E}_\sigma R_n\{f - g : f, g \in \mathcal{F}, P_n(f - g)^2 \leq c_3 r\} \\ &= 2\mathbb{E}_\sigma R_n\{f : f \in \mathcal{F}, P_n f^2 \leq c_3 r/4\}.\end{aligned}$$

Combining with Lemma 6.6 gives

$$2c_1 \mathbb{E}_\sigma R_n\{f \in \mathcal{F} : P_n(f - \hat{f})^2 \leq c_3 r\} + \frac{(c_2 + 2)x}{n}$$

$$\leq 4c_1 \left(\frac{2}{n}\sum_{i=1}^n \min\left\{\frac{c_3 r}{4}, \hat{\lambda}_i\right\}\right)^{1/2} + \frac{(c_2 + 2)x}{n}.$$

Let $\widehat{\psi}_n(r)$ denote the right-hand side. Notice that $\widehat{\psi}_n$ is a sub-root function, so the estimate of Theorem 5.4 can be applied. To compute the fixed point of $B\widehat{\psi}_n$, first notice that adding a constant $a$ to a sub-root function can increase its fixed point by at most $2a$. Thus, it suffices to show that

$$r \leq 4c_1 \left(\frac{2}{n}\sum_{i=1}^n \min\left\{\frac{c_3 r}{4}, \hat{\lambda}_i\right\}\right)^{1/2}$$

implies

(6.2) $$r \leq c \min_{0 \leq h \leq n}\left(\frac{h}{n} + \sqrt{\frac{1}{n}\sum_{i>h}\hat{\lambda}_i}\right)$$

for some universal constant $c$. Under this hypothesis,

$$\begin{aligned}\left(\frac{r}{4c_1}\right)^2 &\leq \frac{2}{n}\sum_{i=1}^n \min\left\{\frac{c_3 r}{4}, \hat{\lambda}_i\right\} \\ &= \frac{2}{n} \min_{S \subseteq \{1,\ldots,n\}}\left(\sum_{i \in S}\frac{c_3 r}{4} + \sum_{i \notin S}\hat{\lambda}_i\right) \\ &= \frac{2}{n} \min_{0 \leq h \leq n}\left(\frac{c_3 h r}{4} + \sum_{i>h}\hat{\lambda}_i\right).\end{aligned}$$

Solving the quadratic inequality for each value of $h$ gives (6.2). □

## APPENDIX

**A.1. Additional material.** This section contains a collection of results that is needed in the proofs. Most of them are classical or easy to derive from classical results. We present proofs for the sake of completeness.

34    P. L. BARTLETT, O. BOUSQUET AND S. MENDELSONRecall the following improvement of Rio's [29] version of Talagrand's concentration inequality, which is due to Bousquet [7, 8].

THEOREM A.1.  *Let $c > 0$, let $X_i$ be independent random variables distributed according to $P$ and let $\mathcal{F}$ be a set of functions from $\mathcal{X}$ to $\mathbb{R}$. Assume that all functions $f$ in $\mathcal{F}$ satisfy $\mathbb{E}f = 0$ and $\|f\|_\infty \leq c$.*

*Let $\sigma$ be a positive real number such that $\sigma^2 \geq \sup_{f \in \mathcal{F}} \mathsf{Var}[f(X_i)]$. Then, for any $x \geq 0$,*

$$\Pr(Z \geq \mathbb{E}Z + x) \leq \exp\left(-vh\left(\frac{x}{cv}\right)\right),$$

*where $Z = \sup_{f \in \mathcal{F}} \sum_{i=1}^n f(X_i)$, $h(x) = (1+x)\log(1+x) - x$ and $v = n\sigma^2 + 2c\mathbb{E}Z$. Also, with probability at least $1 - e^{-x}$,*

$$Z \leq \mathbb{E}Z + \sqrt{2xv} + \frac{cx}{3}.$$

In a similar way one can obtain a concentration result for the Rademacher averages of a class (using the result of [5]; see also [6]). In order to obtain the appropriate constants, notice that

$$\mathbb{E}_\sigma \sup_{f \in \mathcal{F}} \sum_{i=1}^n \sigma_i f(X_i) = \mathbb{E}_\sigma \sup_{f \in \mathcal{F}} \sum_{i=1}^n \sigma_i (f(X_i) - (b-a)/2)$$

and $|f - (b-a)/2| \leq (b-a)/2$.

THEOREM A.2.  *Let $\mathcal{F}$ be a class of functions that map $\mathcal{X}$ into $[a, b]$. Let*

$$Z = \mathbb{E}_\sigma \sup_{f \in \mathcal{F}} \sum_{i=1}^n \sigma_i f(X_i) = n \mathbb{E}_\sigma R_n \mathcal{F}.$$

*Then for all $x \geq 0$,*

$$\Pr\left(Z \geq \mathbb{E}Z + \sqrt{(b-a)x\mathbb{E}Z} + \frac{(b-a)x}{6}\right) \leq e^{-x}$$

*and*

$$\Pr(Z \leq \mathbb{E}Z - \sqrt{(b-a)x\mathbb{E}Z}) \leq e^{-x}.$$

LEMMA A.3.  *For $u, v \geq 0$,*

$$\sqrt{u+v} \leq \sqrt{u} + \sqrt{v},$$

*and for any $\alpha > 0$,*

$$2\sqrt{uv} \leq \alpha u + \frac{v}{\alpha}.$$

LEMMA A.4. *Fix $x > 0$, and let $\mathcal{F}$ be a class of functions with ranges in $[a, b]$. Then, with probability at least $1 - e^{-x}$,*

$$\mathbb{E} R_n \mathcal{F} \leq \inf_{\alpha \in (0,1)} \left( \frac{1}{1-\alpha} \mathbb{E}_\sigma R_n \mathcal{F} + \frac{(b-a)x}{4n\alpha(1-\alpha)} \right).$$

*Also, with probability at least $1 - e^{-x}$,*

$$\mathbb{E}_\sigma R_n \mathcal{F} \leq \inf_{\alpha > 0} \left( (1+\alpha) \mathbb{E} R_n \mathcal{F} + \frac{(b-a)x}{2n} \left( \frac{1}{2\alpha} + \frac{1}{3} \right) \right).$$

PROOF. The second inequality of Theorem A.2 and Lemma A.3 imply that with probability at least $1 - e^{-x}$,

$$\mathbb{E} R_n \mathcal{F} \leq \mathbb{E}_\sigma R_n \mathcal{F} + \sqrt{\frac{(b-a)x}{n} \mathbb{E} R_n \mathcal{F}}$$

$$\leq \mathbb{E}_\sigma R_n \mathcal{F} + \alpha \mathbb{E} R_n \mathcal{F} + \frac{(b-a)x}{4n\alpha},$$

and the first claim of the lemma follows. The proof of the second claim is similar, but uses the first inequality of Theorem A.2. □

A standard fact is that the expected deviation of the empirical means from the actual ones can be controlled by the Rademacher averages of the class.

LEMMA A.5. *For any class of functions $\mathcal{F}$,*

$$\max\left( \mathbb{E} \sup_{f \in \mathcal{F}} (Pf - P_n f), \mathbb{E} \sup_{f \in \mathcal{F}} (P_n f - Pf) \right) \leq 2 \mathbb{E} R_n \mathcal{F}.$$

PROOF. Let $X'_1, \ldots, X'_n$ be an independent copy of $X_1, \ldots, X_n$, and set $P'_n$ to be the empirical measure supported on $X'_1, \ldots, X'_n$. By the convexity of the supremum and by symmetry,

$$\begin{aligned}
\mathbb{E} \sup_{f \in \mathcal{F}} (Pf - P_n f) &= \mathbb{E} \sup_{f \in \mathcal{F}} (\mathbb{E} P'_n f - P_n f) \\
&\leq \mathbb{E} \sup_{f \in \mathcal{F}} (P'_n f - P_n f) \\
&= \frac{1}{n} \mathbb{E} \sup_{f \in \mathcal{F}} \left[ \sum_{i=1}^n \sigma_i f(X'_i) - \sigma_i f(X_i) \right] \\
&\leq \frac{1}{n} \mathbb{E} \sup_{f \in \mathcal{F}} \sum_{i=1}^n \sigma_i f(X'_i) + \frac{1}{n} \mathbb{E} \sup_{f \in \mathcal{F}} \sum_{i=1}^n -\sigma_i f(X_i) \\
&= 2 \mathbb{E} \sup_{f \in \mathcal{F}} R_n f.
\end{aligned}$$

Using an identical argument, the same holds for $P_n f - Pf$. □

In addition, recall the following contraction property, which is due to Ledoux and Talagrand [17].

THEOREM A.6. *Let $\phi$ be a contraction, that is, $|\phi(x) - \phi(y)| \leq |x - y|$. Then, for every class $\mathcal{F}$,*

$$\mathbb{E}_\sigma R_n \phi \circ \mathcal{F} \leq \mathbb{E}_\sigma R_n \mathcal{F},$$

*where $\phi \circ \mathcal{F} := \{\phi \circ f : f \in \mathcal{F}\}$.*

The interested reader may find some additional useful properties of the Rademacher averages in [3, 27].

### A.2. Proofs.

PROOF OF THEOREM 2.1. Define $V^+ = \sup_{f \in \mathcal{F}}(Pf - P_n f)$. Since $\sup_{f \in \mathcal{F}} \mathsf{Var}[f(X_i)] \leq r$, and $\|f - Pf\|_\infty \leq b - a$, Theorem A.1 implies that, with probability at least $1 - e^{-x}$,

$$V^+ \leq \mathbb{E}V^+ + \sqrt{\frac{2xr}{n} + \frac{4x(b-a)\mathbb{E}V^+}{n}} + \frac{(b-a)x}{3n}.$$

Thus by Lemma A.3, with probability at least $1 - e^{-x}$,

$$V^+ \leq \inf_{\alpha > 0}\left((1+\alpha)\mathbb{E}V^+ + \sqrt{\frac{2rx}{n}} + (b-a)\left(\frac{1}{3} + \frac{1}{\alpha}\right)\frac{x}{n}\right).$$

Applying Lemma A.5 gives the first assertion of Theorem 2.1. The second part of the theorem follows by combining the first one and Lemma A.4, and noticing that $\inf_\alpha f(\alpha) + \inf_\alpha g(\alpha) \leq \inf_\alpha (f(\alpha) + g(\alpha))$. Finally, the fact that the same results hold for $\sup_{f \in \mathcal{F}}(P_n f - Pf)$ can be easily obtained by applying the above reasoning to the class $-\mathcal{F} = \{-f : f \in \mathcal{F}\}$ and noticing that the Rademacher averages of $-\mathcal{F}$ and $\mathcal{F}$ are identical. □

PROOF OF LEMMA 3.2. To prove the continuity of $\psi$, let $x > y > 0$, and note that since $\psi$ is nondecreasing, $|\psi(x) - \psi(y)| = \psi(x) - \psi(y)$. From the fact that $\psi(r)/\sqrt{r}$ is nonincreasing it follows that $\psi(x)/\sqrt{y} \leq \sqrt{x}\psi(y)/y$, and thus

$$\psi(x) - \psi(y) = \sqrt{y}\frac{\psi(x)}{\sqrt{y}} - \psi(y) \leq \psi(y)\frac{\sqrt{x} - \sqrt{y}}{\sqrt{y}}.$$

Letting $x$ tend to $y$, $|\psi(x) - \psi(y)|$ tends to 0, and $\psi$ is left-continuous at $y$. A similar argument shows the right-sided continuity of $\psi$.

As for the second part of the claim, note that $\psi(x)/x$ is nonnegative and continuous on $(0, \infty)$, and since $1/\sqrt{x}$ is strictly decreasing on $(0, \infty)$, then $\psi(x)/x$ is also strictly decreasing.

Observe that if $\psi(x)/x$ is always larger than 1 on $(0, \infty)$, then $\lim_{x \to \infty} \psi(x)/\sqrt{x} = \infty$, which is impossible. On the other hand, if $\psi(x)/x < 1$ on $(0, \infty)$, then $\lim_{x \to 0} \psi(x)/\sqrt{x} = 0$, contrary to the assumption that $\psi$ is nontrivial. Thus the equation $\psi(r)/r = 1$ has a positive solution and this solution is unique by monotonicity.

Finally, if for some $r > 0$, $r \geq \psi(r)$, then $\psi(t)/t \leq 1$ for all $t \geq r$ [since $\psi(x)/x$ is nonincreasing] and thus $r^* \leq r$. The other direction follows in a similar manner. $\square$

PROOF OF LEMMA 3.4. Observe that, by symmetry of the Rademacher random variables, one has $\psi(r) = \mathbb{E}_\sigma R_n \{f - \hat{f} : f \in \mathcal{F}, T(f - \hat{f}) \leq r\}$ so that, by translating the class, it suffices to consider the case where $\hat{f} = 0$.

Note that $\psi$ is nonnegative, since by Jensen's inequality

$$\mathbb{E}_\sigma \sup_{f \in \mathcal{F}} R_n f \geq \sup_{f \in \mathcal{F}} \mathbb{E}_\sigma R_n f = 0.$$

Moreover, $\psi$ is nondecreasing since $\{f \in \mathcal{F} : T(f) \leq r\} \subset \{f \in \mathcal{F} : T(f) \leq r'\}$ for $r \leq r'$. It remains to show that for any $0 < r_1 \leq r_2$, $\psi(r_1) \geq \sqrt{r_1/r_2} \cdot \psi(r_2)$. To this end, fix any sample and any realization of the Rademacher random variables, and set $f_0$ to be a function for which

$$\sup_{f \in \mathcal{F}, T(f) \leq r_2} \sum_{i=1}^{n} \sigma_i f(x_i)$$

is attained (if the supremum is not attained only a slight modification is required). Since $T(f_0) \leq r_2$, then $T(\sqrt{r_1/r_2} \cdot f_0) \leq r_1$ by assumption. Furthermore, since $\mathcal{F}$ is star-shaped, the function $\sqrt{r_1/r_2} f_0$ belongs to $\mathcal{F}$ and satisfies that $T(\sqrt{r_1/r_2} f_0) \leq r_1$. Hence

$$\sup_{f \in \mathcal{F} : T(f) \leq r_1} \sum_{i=1}^{n} \sigma_i f(x_i) \geq \sum_{i=1}^{n} \sigma_i \sqrt{\frac{r_1}{r_2}} \cdot f_0(x_i)$$

$$= \sqrt{\frac{r_1}{r_2}} \sup_{f \in \mathcal{F} : T(f) \leq r_2} \sum_{i=1}^{n} \sigma_i f(x_i),$$

and the result follows by taking expectations with respect to the Rademacher random variables. $\square$

PROOF OF COROLLARY 3.7. The proof uses the following result of [11], which relates the empirical Rademacher averages to the empirical $L_2$ entropy of the class. The covering number $\mathcal{N}(\epsilon, \mathcal{F}, L_2(P_n))$ is the cardinality of the smallest subset $\hat{\mathcal{F}}$ of $L_2(P_n)$ for which every element of $\mathcal{F}$ is within $\epsilon$ of some element of $\hat{\mathcal{F}}$.

THEOREM B.7 ([11]). *There exists an absolute constant $C$ such that for every class $\mathcal{F}$ and every $X_1, \ldots, X_n \in \mathcal{X}$,*

$$\mathbb{E}_\sigma R_n \mathcal{F} \leq \frac{C}{\sqrt{n}} \int_0^\infty \sqrt{\log \mathcal{N}(\varepsilon, \mathcal{F}, L_2(P_n))} \, d\varepsilon.$$

Define the sub-root function

$$\psi(r) = 10 \mathbb{E} R_n \{f \in \mathrm{star}(\mathcal{F}, 0) : Pf^2 \leq r\} + \frac{11 \log n}{n}.$$

If $r \geq \psi(r)$, then Corollary 2.2 implies that, with probability at least $1 - 1/n$,

$$\{f \in \mathrm{star}(\mathcal{F}, 0) : Pf^2 \leq r\} \subseteq \{f \in \mathrm{star}(\mathcal{F}, 0) : P_n f^2 \leq 2r\},$$

and thus

$$\mathbb{E} R_n \{f \in \mathrm{star}(\mathcal{F}, 0) : Pf^2 \leq r\} \leq \mathbb{E} R_n \{f \in \mathrm{star}(\mathcal{F}, 0) : P_n f^2 \leq 2r\} + \frac{1}{n}.$$

It follows that $r^* = \psi(r^*)$ satisfies

$$\text{(A.1)} \qquad r^* \leq 10 \mathbb{E} R_n \{f \in \mathrm{star}(\mathcal{F}, 0) : P_n f^2 \leq 2r^*\} + \frac{1 + 11 \log n}{n}.$$

But Theorem B.7 shows that

$$\mathbb{E} R_n \{f \in \mathrm{star}(\mathcal{F}, 0) : P_n f^2 \leq 2r^*\}$$

$$\leq \frac{C}{\sqrt{n}} \mathbb{E} \int_0^{\sqrt{2r^*}} \sqrt{\log \mathcal{N}(\varepsilon, \mathrm{star}(\mathcal{F}, 0), L_2(P_n))} \, d\varepsilon.$$

It is easy to see that we can construct an $\epsilon$-cover for $\mathrm{star}(\mathcal{F}, 0)$ using an $\epsilon/2$-cover for $\mathcal{F}$ and an $\epsilon/2$-cover for the interval $[0, 1]$, which implies

$$\log \mathcal{N}(\varepsilon, \mathrm{star}(\mathcal{F}, 0), L_2(P_n)) \leq \log \mathcal{N}\left(\frac{\varepsilon}{2}, \mathcal{F}, L_2(P_n)\right)\left(\left\lceil \frac{2}{\epsilon} \right\rceil + 1\right).$$

Now, recall that [14] for any probability distribution $P$ and any class $\mathcal{F}$ with VC-dimension $d < \infty$,

$$\log \mathcal{N}\left(\frac{\varepsilon}{2}, \mathcal{F}, L_2(P)\right) \leq cd \log\left(\frac{1}{\epsilon}\right).$$

Therefore

$$\mathbb{E} R_n \{f \in \mathrm{star}(\mathcal{F}, 0) : P_n f^2 \leq 2r^*\} \leq \sqrt{\frac{cd}{n}} \int_0^{\sqrt{2r^*}} \sqrt{\log\left(\frac{1}{\epsilon}\right)} \, d\varepsilon$$

$$\leq \sqrt{\frac{cd r^* \log(1/r^*)}{n}}$$

$$\leq \sqrt{c\left(\frac{d^2}{n^2} + \frac{dr^* \log(n/ed)}{n}\right)},$$

where $c$ represents an absolute constant whose value may change from line to line. Substituting into (A.1) and solving for $r^*$ shows that

$$r^* \leq \frac{cd\log(n/d)}{n},$$

provided $n \geq d$. The result follows from Theorem 3.3. □

PROOF OF THEOREM 5.2. Let $f^* = \arg\min_{f \in \mathcal{F}} P\ell_f$. (For simplicity, assume that the minimum exists; if it does not, the proof is easily extended by considering the limit of a sequence of functions with expected loss approaching the infimum.) Then, by definition of $\hat{f}$, $P_n\ell_{\hat{f}} \leq P_n\ell_{f^*}$. Since the variance of $\ell_{f^*}(X_i, Y_i)$ is no more than some constant times $L^*$, we can apply Bernstein's inequality (see, e.g., [10], Theorem 8.2) to show that with probability at least $1 - e^{-x}$,

$$P_n\ell_{\hat{f}} \leq P_n\ell_{f^*} \leq P\ell_{f^*} + c\left(\sqrt{\frac{P\ell_{f^*}x}{n}} + \frac{x}{n}\right) = L^* + c\left(\sqrt{\frac{L^*x}{n}} + \frac{x}{n}\right).$$

Thus, by Theorem 3.3, with probability at least $1 - 2e^{-x}$,

$$P\ell_{\hat{f}} \leq \frac{K}{K-1}\left(L^* + c\left(\sqrt{\frac{L^*x}{n}} + \frac{x}{n}\right)\right) + cK\left(r^* + \frac{x}{n}\right).$$

Setting

$$K - 1 = \sqrt{\frac{\max(L^*, x/n)}{r^*}},$$

noting that $r^* \geq x/n$ and simplifying gives the first inequality. A similar argument using Theorem 4.1 implies the second inequality. □

PROOF OF LEMMA 6.6. Introduce the operator $\hat{C}_n$ on $\mathcal{H}$ defined by

$$(\hat{C}_n f)(x) = \frac{1}{n}\sum_{i=1}^n f(X_i)k(X_i, x),$$

so that, using (6.1),

$$\langle g, \hat{C}_n f\rangle = \frac{1}{n}\sum_{i=1}^n f(X_i)g(X_i),$$

and $\langle f, \hat{C}_n f\rangle = P_n f^2$, implying that $\hat{C}_n$ is positive semidefinite.

Suppose that $f$ is an eigenfunction of $\hat{C}_n$ with eigenvalue $\lambda$. Then for all $i$

$$\lambda f(X_i) = (\hat{C}_n f)(X_i) = \frac{1}{n}\sum_{j=1}^n f(X_j)k(X_j, X_i).$$

Thus, the vector $(f(X_1), \ldots, f(X_n))$ is either zero (which implies $\hat{C}_n f = 0$ and hence $\lambda = 0$) or is an eigenvector of $\hat{T}_n$ with eigenvalue $\lambda$. Conversely, if $\hat{T}_n v = \lambda v$ for some vector $v$, then

$$\hat{C}_n \left( \sum_{i=1}^n v_i k(X_i, \cdot) \right) = \frac{1}{n} \sum_{i,j=1}^n v_i k(X_i, X_j) k(X_j, \cdot) = \frac{\lambda}{n} \sum_{j=1}^n v_j k(X_j, \cdot).$$

Thus, the eigenvalues of $\hat{T}_n$ are the same as the $n$ largest eigenvalues of $\hat{C}_n$, and the remaining eigenvalues of $\hat{C}_n$ are zero. Let $(\hat{\lambda}_i)$ denote these eigenvalues, arranged in a nonincreasing order.

Let $(\Phi_i)_{i \geq 1}$ be an orthonormal basis of $\mathcal{H}$ of eigenfunctions of $\hat{C}_n$ (such that $\Phi_i$ is associated with $\hat{\lambda}_i$). Fix $0 \leq h \leq n$ and note that for any $f \in \mathcal{H}$

$$\sum_{i=1}^n \sigma_i f(X_i) = \left\langle f, \sum_{i=1}^n \sigma_i k(X_i, \cdot) \right\rangle$$

$$= \left\langle \sum_{j=1}^h \sqrt{\hat{\lambda}_j} \langle f, \Phi_j \rangle \Phi_j, \sum_{j=1}^h \frac{1}{\sqrt{\hat{\lambda}_j}} \left\langle \sum_{i=1}^n \sigma_i k(X_i, \cdot), \Phi_j \right\rangle \Phi_j \right\rangle$$

$$+ \left\langle f, \sum_{j>h} \left\langle \sum_{i=1}^n \sigma_i k(X_i, \cdot), \Phi_j \right\rangle \Phi_j \right\rangle.$$

If $\|f\| \leq 1$ and

$$r \geq P_n f^2 = \langle f, \hat{C}_n f \rangle = \sum_{i \geq 1} \hat{\lambda}_i \langle f, \Phi_i \rangle^2,$$

then by the Cauchy–Schwarz inequality

$$\sum_{i=1}^n \sigma_i f(X_i) \leq \sqrt{r \sum_{j=1}^h \frac{1}{\hat{\lambda}_j} \left\langle \sum_{i=1}^n \sigma_i k(X_i, \cdot), \Phi_j \right\rangle^2}$$

(A.2)

$$+ \sqrt{\sum_{j>h} \left\langle \sum_{i=1}^n \sigma_i k(X_i, \cdot), \Phi_j \right\rangle^2}.$$

Moreover,

$$\frac{1}{n} \mathbb{E}_\sigma \left\langle \sum_{i=1}^n \sigma_i k(X_i, \cdot), \Phi_j \right\rangle^2 = \frac{1}{n} \mathbb{E}_\sigma \sum_{i,\ell=1}^n \sigma_i \sigma_\ell \langle k(X_i, \cdot), \Phi_j \rangle \langle k(X_l, \cdot), \Phi_j \rangle$$

$$= \frac{1}{n} \sum_{i=1}^n \langle k(X_i, \cdot), \Phi_j \rangle^2$$

$$= \langle \Phi_j, \hat{C}_n \Phi_j \rangle$$
$$= \hat{\lambda}_j.$$

Using (A.2) and Jensen's inequality, it follows that

$$\mathbb{E}_\sigma R_n\{f \in \mathcal{F} : P_n f^2 \leq r\} \leq \frac{1}{\sqrt{n}} \min_{0 \leq h \leq n} \left\{ \sqrt{hr} + \sqrt{\sum_{j=h+1}^n \hat{\lambda}_j} \right\},$$

which implies the result. $\square$

**Acknowledgments.** We are grateful to the anonymous reviewers and to Vu Ha for comments that improved the paper. We are very much indebted to Regis Vert for suggestions that led to an important improvement of the results and a simplification of the proofs.

P. L. BARTLETT  
DEPARTMENT OF STATISTICS AND  
  DIVISION OF COMPUTER SCIENCE  
367 EVANS HALL  
UNIVERSITY OF CALIFORNIA AT BERKELEY  
BERKELEY, CALIFORNIA 94720-3860  
USA  
E-MAIL: bartlett@stat.berkeley.edu  

O. BOUSQUET  
EMPIRICAL INFERENCE DEPARTMENT  
MAX PLANCK INSTITUTE FOR  
  BIOLOGICAL CYBERNETICS  
SPEMANNSTR. 38  
D-72076 TÜBINGEN  
GERMANY  
E-MAIL: olivier.bousquet@tuebingen.mpg.de  

S. MENDELSON  
CENTRE FOR MATHEMATHICS  
  AND ITS APPLICATIONS  
INSTITUTE OF ADVANCED STUDIES  
AUSTRALIAN NATIONAL UNIVERSITY  
CANBERRA, ACT 0200  
AUSTRALIA  
E-MAIL: shahar.mendelson@anu.edu.au